# Basis of Representation of Universal Algebra

## Aleks Kleyn

Abstract. We say that there is a representation of the universal algebra $B$ in the universal algebra $A$ if the set of endomorphisms of the universal algebra $A$ has the structure of universal algebra $B$. Therefore, the role of representation of the universal algebra is similar to the role of symmetry in geometry and physics. Morphism of the representation is the mapping that conserves the structure of the representation. Exploring of morphisms of the representation leads to the concepts of generating set and basis of representation. The set of automorphisms of the representation of the universal algebra forms the group. Twin representations of this group in basis manifold of the representation are called active and passive representations. Passive representation in basis manifold is underlying of concept of geometric object and the theory of invariants of the representation of the universal algebra.

## Contents

## 1. Preface

This paper is based on the chapter [3]-3.

The role of representation of the universal algebra is similar to the role of symmetry in geometry and physics. In both cases we study the structured set of transformations; knowledge of the structure of this set gives us additional information about the object under study.

Aleks_Kleyn@MailAPS.org.
http://sites.google.com/site/AleksKleyn/.
http://arxiv.org/a/kleyn_a_1.
http://AleksKleyn.blogspot.com/.

We consider the theory of representations of universal algebra as an extension of the theory of universal algebra. Any algebraic structure assumes a set of mappings preserving this structure. Mapping that preserves the structure of the representation of universal algebra, is called a morphism of representations.

Automorphism of representation is special case of morphism of representations. The study of automorphisms representation is directly related to the necessity to answer the question what is the structure of the generating set of representation. If the representation has minimal generating set, then such set is called basis of the representation. Automorphism of representation maps basis into basis.

The set of automorphisms of representation forms a group. Twin representations of this group in basis manifold of the representation are called active and passive representations. Passive representation in basis manifold is underlying of concept of geometric object and the theory of invariants of the representation of the universal algebra.

## 2. Conventions

**Convention 2.1.** *In [4], an arbitrary operation of algebra is denoted by letter $\omega$, and $\Omega$ is the set of operations of some universal algebra. Correspondingly, the universal algebra with the set of operations $\Omega$ is denoted as $\Omega$-algebra. Similar notations we see in [2] with small difference that an operation in the algebra is denoted by letter $f$ and $\mathcal{F}$ is the set of operations. I preferred first case of notations because in this case it is easier to see where I use operation.*                □

**Convention 2.2.** *Let $A$ be $\Omega_1$-algebra. Let $B$ be $\Omega_2$-algebra. Notation*

$$A \longrightarrow_* B$$

*means that there is representation of $\Omega_1$-algebra $A$ in $\Omega_2$-algebra $B$.*                □

Without a doubt, the reader may have questions, comments, objections. I will appreciate any response.

## 3. Generating Set of Representation

**Definition 3.1.** Let

$$f : A \to {}^*M$$

be representation of $\Omega_1$-algebra $A$ in $\Omega_2$-algebra $M$. The set $N \subset M$ is called **stable set of representation** $f$, if $f(a)(m) \in N$ for each $a \in A$, $m \in N$.                □

We also say that the set $M$ is stable with respect to the representation $f$.

**Theorem 3.2.** *Let*

$$f : A \to {}^*M$$

*be representation of $\Omega_1$-algebra $A$ in $\Omega_2$-algebra $M$. Let set $N \subset M$ be subalgebra of $\Omega_2$-algebra $M$ and stable set of representation $f$. Then there exists representation*

$$f_N : A \to {}^*N$$

*such that $f_N(a) = f(a)|_N$. Representation $f_N$ is called **subrepresentation of representation** $f$.*

*Proof.* Let $\omega_1$ be $n$-ary operation of $\Omega_1$-algebra $A$. Then for each $a_1, ..., a_n \in A$ and each $b \in N$

$$(f_N(a_1)...f_N(a_n)\omega_1)(b) = (f(a_1)...f(a_n)\omega_1)(b)$$
$$= f(a_1...a_n\omega_1)(b)$$
$$= f_N(a_1...a_n\omega_1)(b)$$

Let $\omega_2$ be $n$-ary operation of $\Omega_2$-algebra $M$. Then for each $b_1, ..., b_n \in N$ and each $a \in A$

$$f_N(a)(b_1)...f_N(a)(b_n)\omega_2 = f(a)(b_1)...f(a)(b_n)\omega_2$$
$$= f(a)(b_1...b_n\omega_2)$$
$$= f_N(a)(b_1...b_n\omega_2)$$

We proved the statement of theorem. $\square$

From the theorem 3.2, it follows that if $f_N$ is subrepresentation of representation $f$, then the mapping $(id : A \to A, id_N : N \to M)$ is morphism of representations.

**Theorem 3.3.** *The set[1] $\mathcal{B}_f$ of all subrepresentations of representation $f$ generates a closure system on $\Omega_2$-algebra $M$ and therefore is a complete lattice.*

*Proof.* Let $(K_\lambda)_{\lambda \in \Lambda}$ be the set off subalgebras of $\Omega_2$-algebra $M$ that are stable with respect to representation $f$. We define the operation of intersection on the set $\mathcal{B}_f$ according to rule

$$\bigcap f_{K_\lambda} = f_{\cap K_\lambda}$$

We defined the operation of intersection of subrepresentations properly. $\cap K_\lambda$ is subalgebra of $\Omega_2$-algebra $M$. Let $m \in \cap K_\lambda$. For each $\lambda \in \Lambda$ and for each $a \in A$, $f(a)(m) \in K_\lambda$. Therefore, $f(a)(m) \in \cap K_\lambda$. Therefore, $\cap K_\lambda$ is the stable set of representation $f$. $\square$

We denote the corresponding closure operator by $J(f)$. Thus $J(f, X)$ is the intersection of all subalgebras of $\Omega_2$-algebra $M$ containing $X$ and stable with respect to representation $f$.

**Theorem 3.4.** *Let[2]*

$$f : A \to {}^*M$$

*be representation of $\Omega_1$-algebra $A$ in $\Omega_2$-algebra $M$. Let $X \subset M$. Define a subset $X_k \subset M$ by induction on $k$.*

$$X_0 = X$$
$$x \in X_k => x \in X_{k+1}$$
$$x_1 \in X_k, ..., x_n \in X_k, \omega \in \Omega_2(n) => x_1...x_n\omega \in X_{k+1}$$
$$x \in X_k, a \in A => f(a)(x) \in X_{k+1}$$

*Then*

$$\bigcup_{k=0}^{\infty} X_k = J(f, X)$$

---

[1] This definition is similar to definition of the lattice of subalgebras ([4], p. 79, 80)

[2] The statement of theorem is similar to the statement of theorem 5.1, [4], p. 79.

*Proof.* If we put $U = \cup X_k$, then by definition of $X_k$, we have $X_0 \subset J(f, X)$, and if $X_k \subset J(f, X)$, then $X_{k+1} \subset J(f, X)$. By induction it follows that $X_k \subset J(f, X)$ for all $k$. Therefore,

$$(3.1) \qquad\qquad U \subset J(f, X)$$

If $a \in U^n$, $a = (a_1, ..., a_n)$, where $a_i \in X_{k_i}$, and if $k = \max\{k_1, ..., k_n\}$, then $a_1...a_n\omega \in X_{k+1} \subset U$. Therefore, $U$ is subalgebra of $\Omega_2$-algebra $M$.

If $m \in U$, then there exists such $k$ that $m \in X_k$. Therefore, $f(a)(m) \in X_{k+1} \subset U$ for any $a \in A$. Therefore, $U$ is stable set of the representation $f$.

Since $U$ is subalgebra of $\Omega_2$-algebra $M$ and is a stable set of the representation $f$, then subrepresentation $f_U$ is defined. Therefore,

$$(3.2) \qquad\qquad J(f, X) \subset U$$

From (3.1), (3.2), it follows that $J(f, X) = U$.                         $\square$

**Definition 3.5.** $J(f, X)$ is called **subrepresentation generated by set** $X$, and $X$ is a **generating set of subrepresentation** $J(f, X)$. In particular, a **generating set of representation** $f$ is a subset $X \subset M$ such that $J(f, X) = M$.                         $\square$

It is easy to see that the definition of generating set of representation does not depend on whether representation is effective or not. For this reason hereinafter we will assume that the representation is effective and we will use convention for effective $T\star$-representation in remark [3]-2.1.9. We also will use notation

$$R \circ m = R(m)$$

for image of $m \in M$ under the endomorphism of effective representation. According to the definition of product of mappings, for any endomorphisms $R$, $S$ the following equation is true

$$(3.3) \qquad\qquad (R \circ S) \circ m = R \circ (S \circ m)$$

The equation (3.3) is associative law for $\circ$ and allows us to write expression

$$R \circ S \circ m$$

without brackets.

From theorem 3.4, it follows next definition.

**Definition 3.6.** Let $X \subset M$. For each $x \in J(f, X)$ there exists $\Omega_2$-word defined according to following rules.

(1) If $m \in X$, then $m$ is $\Omega_2$-word.
(2) If $m_1, ..., m_n$ are $\Omega_2$-words and $\omega \in \Omega_2(n)$, then $m_1...m_n\omega$ is $\Omega_2$-word.
(3) If $m$ is $\Omega_2$-word and $a \in A$, then $am$ is $\Omega_2$-word.

$\Omega_2$-**word** $w(f, X, m)$ represents given element $m \in J(f, X)$. We will identify an element $m \in J(f, X)$ and corresponding it $\Omega_2$-word using equation

$$m = w(f, X, m)$$

Similarly, for an arbitrary set $B \subset J(f, X)$ we consider the set of $\Omega_2$-words[3]

$$w(f, X, B) = \{w(f, X, m) : m \in B\}$$

---

[3]The expression $w(f, X, m)$ is a special case of the expression $w(f, X, B)$, namely
$$w(f, X, \{m\}) = \{w(f, X, m)\}$$

We also use notation

$$w(f, X, B) = (w(f, X, m), m \in B)$$

Denote $w(f, X)$ **the set of $\Omega_2$-words of representation** $J(f, X)$. □

**Theorem 3.7.** *Endomorpism $R$ of representation*

$$f : A \to {}^*M$$

*generates the mapping of $\Omega_2$-words*

$$w[f, X, R] : w(f, X) \to w(f, X') \quad X \subset M \quad X' = R \circ X$$

*such that*

(1) *If $m \in X$, $m' = R \circ m$, then*

$$w[f, X, R](m) = m'$$

(2) *If*

$$m_1, ..., m_n \in w(f, X)$$

$$m'_1 = w[f, X, R](m_1) \quad ... \quad m'_n = w[f, X, R](m_n)$$

*then for operation $\omega \in \Omega_2(n)$ holds*

$$w[f, X, R](m_1...m_n\omega) = m'_1...m'_n\omega$$

(3) *If*

$$m \in w(f, X) \quad m' = w[f, X, R](m) \quad a \in A$$

*then*

$$w[f, X, R](am) = am'$$

*Proof.* Statements (1), (2) of the theorem are true by definition of the endomorpism $R$. Because $r = \mathrm{id}$, the statement (3) of the theorem follows from the equation [3]-(2.2.4). □

*Remark* 3.8. Let $R$ be endomorphism of representation $f$. Let

$$m \in J(f, X) \quad m' = R \circ m \quad X' = R \circ X$$

The theorem 3.7 states that $m' \in J_f(X')$ . The theorem 3.7 also states that $\Omega_2$-word representing $m$ relative $X$ and $\Omega_2$-word representing $m'$ relative $X'$ are generated according to the same algorithm. This allows considering of the set of $\Omega_2$-words $w(f, X', m')$ as mapping

$$W(f, X, m) : X' \to w(f, X', m')$$

$$W(f, X, m)(X') = W(f, X, m) \circ X'$$

such that, if for certain endomorphism $R$

$$X' = R \circ X \quad m' = R \circ m$$

then

$$W(f, X, m) \circ X' = w(f, X', m') = m'$$

The mapping $W(f, X, m)$ is called **coordinates of element $m$ relative to set $X$**. Similarly, we consider coordinates of a set $B \subset J(f, X)$ relative to the set $X$

$$W(f, X, B) = \{W(f, X, m) : m \in B\} = (W(f, X, m), m \in B)$$

Denote $W(f, X)$ **the set of coordinates of representation** $J(f, X)$. $\qquad \square$

**Theorem 3.9.** *There is a structure of $\Omega_2$-algebra on the set of coordinates $W(f, X)$.*

*Proof.* Let $\omega \in \Omega_2(n)$. Then for any $m_1, ..., m_n \in J(f, X)$, we assume

$$(3.4) \qquad W(f, X, m_1)...W(f, X, m_n)\omega = W(f, X, m_1...m_n\omega)$$

According to the remark 3.8,

$$(3.5) \qquad \begin{aligned} (W(f, X, m_1)...W(f, X, m_n)\omega) \circ X &= W(f, X, m_1...m_n\omega) \circ X \\ &= w(f, X, m_1...m_n\omega) \end{aligned}$$

follows from the equation (3.4). According to rule (2) of the definition 3.6, from the equation (3.5), it follows that

(3.6)
$$\begin{aligned} (W(f, X, m_1)...W(f, X, m_n)\omega) \circ X &= w(f, X, m_1)...w(f, X, m_n)\omega \\ &= (W(f, X, m_1) \circ X)...(W(f, X, m_n) \circ X)\omega \end{aligned}$$

From the equation (3.6), it follows that the operation $\omega$ defined by the equation (3.4) on the set of coordinates is defined properly. $\qquad \square$

**Theorem 3.10.** *There exists the representation of $\Omega_1$-algebra $A$ in $\Omega_2$-algebra $W(f, X)$.*

*Proof.* Let $a \in A$. Then for any $m \in J(f, X)$ we assume

$$(3.7) \qquad aW(f, X, m) = W(f, X, am)$$

According to the remark 3.8,

$$(3.8) \qquad (aW(f, X, m)) \circ X = W(f, X, am) \circ X = w(f, X, am)$$

follows from the equation (3.7). According to rule (3) of the definition 3.6, from the equation (3.8), it follows that

$$(3.9) \qquad (aW(f, X, m)) \circ X = aw(f, X, m) = a(W(f, X, m) \circ X)$$

From the equation (3.9), it follows that the representation (3.7) of $\Omega_1$-algebra $A$ in $\Omega_2$-algebra $W(f, X)$ is defined properly. $\qquad \square$

**Theorem 3.11.** *Let*

$$f : A \to {}^*M$$

*be representation of $\Omega_1$-algebra $A$ in $\Omega_2$-algebra $M$. For given sets $X \subset M$, $X' \subset M$, let map*

$$R_1 : X \to X'$$

*agree with the structure of representation $f$, i. e.*

$$\omega \in \Omega_2(n) \ x_1, \ ..., \ x_n, \ x_1...x_n\omega \in X, \ R_1(x_1...x_n\omega) \in X'$$

$$=> R_1(x_1...x_n\omega) = R_1(x_1)...R_1(x_n)\omega$$

$$x \in X, \ a \in A, \ R_1(ax) \in X'$$

$$=> R_1(ax) = aR_1(x)$$

*Consider the mapping of $\Omega_2$-words*

$$w[f, X, R_1] : w(f, X) \to w(f, X')$$

*that satisfies conditions* (1), (2), (3) *of the theorem 3.7 and such that*

$$x \in X => w[f, X, R_1](x) = R_1(x)$$

*There exists unique endomorphism of $\Omega_2$-algebra $M$*

$$R : M \to M$$

*defined by rule*

$$R \circ m = w[f, X, R_1](w(f, X, m))$$

*which is the morphism of representations $J(f, X)$ and $J(f, X')$.*

*Proof.* We prove the theorem by induction over complexity of $\Omega_2$-word.

If $w(f, X, m) = m$, then $m \in X$. According to condition (1) of theorem 3.7,

$$R \circ m = w[f, X, R_1](w(f, X, m)) = w[f, X, R_1](m) = R_1(m)$$

Therefore, mappings $R$ and $R_1$ coincide on the set $X$, and the mapping $R$ agrees with structure of representation $f$.

Let $\omega \in \Omega_2(n)$. Let the mapping $R$ be defined for $m_1, ..., m_n \in J(f, X)$. Let

$$w_1 = w(f, X, m_1) \quad ... \quad w_n = w(f, X, m_n)$$

If $m = m_1...m_n\omega$, then according to rule (2) of definition 3.6,

$$w(f, X, m) = w_1...w_n\omega$$

According to condition (2) of theorem 3.7,

$$R \circ m = w[f, X, R_1](w(f, X, m)) = w[f, X, R_1](w_1...w_n\omega)$$

$$= w[f, X, R_1](w_1)...w[f, X, R_1](w_n)\omega$$

$$= (R \circ m_1)...(R \circ m_n)\omega$$

Therefore, the mapping $R$ is endomorphism of $\Omega_2$-algebra $M$.

Let the mapping $R$ be defined for $m_1 \in J(f, X)$, $w_1 = w(f, X, m_1)$. Let $a \in A$. If $m = am_1$, then according to rule (3) of definition 3.6,

$$w(f, X, am_1) = aw_1$$

According to condition (3) of theorem 3.7,

$$R \circ m = w[f, X, R_1](w(f, X, m)) = w[f, X, R_1](aw_1)$$

$$= aw[f, X, R_1](w_1) = aR \circ m_1$$

From equation [3]-(2.2.4), it follows that the mapping $R$ is morphism of the representation $f$.

The statement that the endomorphism $R$ is unique and therefore this endomorphism is defined properly follows from the following argument. Let $m \in M$ have different $\Omega_2$-words relative the set $X$, for instance

$$(3.10) \qquad m = x_1 \ldots x_n \omega = ax$$

Because $R$ is endomorphism of representation, then, from the equation (3.10), it follows that

$$(3.11) \qquad R \circ m = R \circ (x_1 \ldots x_n \omega) = (R \circ x_1) \ldots (R \circ x_n)\omega = R \circ (ax) = aR \circ x$$

From the equation (3.11), it follows that

$$(3.12) \qquad R \circ m = (R \circ x_1) \ldots (R \circ x_n)\omega = aR \circ x$$

From equations (3.10), (3.12), it follows that the equation (3.10) is preserved under the mapping. Therefore, the image of $m$ does not depend on the choice of coordinates. $\qquad\square$

*Remark* 3.12. The theorem 3.11 is the theorem of extension of mapping. The only statement we know about the set $X$ is the statement that $X$ is generating set of the representation $f$. However, between the elements of the set $X$ there may be relationships generated by either operations of $\Omega_2$-algebra $M$, or by transformation of representation $f$. Therefore, any mapping of set $X$, in general, cannot be extended to an endomorphism of representation $f$.[4] However, if the mapping $R_1$ is coordinated with the structure of representation on the set $X$, then we can construct an extension of this mapping and this extension is endomorphism of representation $f$. $\qquad\square$

**Definition 3.13.** Let $X$ be the generating set of the representation $f$. Let $R$ be the endomorphism of the representation $f$. The set of coordinates $W(f, X, R \circ X)$ is called **coordinates of endomorphism of representation**. $\qquad\square$

**Definition 3.14.** Let $X$ be the generating set of the representation $f$. Let $R$ be the endomorphism of the representation $f$. Let $m \in M$. We define **superposition of coordinates** of the representation $f$ and the element $m$ as coordinates defined according to rule

$$(3.13) \qquad W(f, X, m) \circ W(f, X, R \circ X) = W(f, X, R \circ m)$$

Let $Y \subset M$. We define superposition of coordinates of the representation $f$ and the set $Y$ according to rule

$$(3.14) \qquad W(f, X, Y) \circ W(f, X, R \circ X) = (W(f, X, m) \circ W(f, X, R \circ X), m \in Y)$$

$$\square$$

**Theorem 3.15.** *Endomorphism $R$ of representation*

$$f : A \to {}^*M$$

*generates the mapping of coordinates of representation*

$$(3.15) \qquad W[f, X, R] : W(f, X) \to W(f, X)$$

---

[4]In the theorem 4.7, requirements to generating set are more stringent. Therefore, the theorem 4.7 says about extension of arbitrary mapping. A more detailed analysis is given in the remark 4.9.

*such that*

$$(3.16) \quad W(f, X, m) \to W[f, X, R] \star W(f, X, m) = W(f, X, R \circ m)$$

$$W[f, X, R] \star W(f, X, m) = W(f, X, m) \circ W(f, X, R \circ X)$$

*Proof.* According to the remark 3.8, we consider equations (3.13), (3.15) relative to given generating set $X$. The word

$$(3.17) \quad W(f, X, m) \circ X = w(f, X, m)$$

corresponds to coordinates $W(f, X, m)$; the word

$$(3.18) \quad W(f, X, R \circ m) \circ X = w(f, X, R \circ m)$$

corresponds to coordinates $W(f, X, R \circ m)$. Therefore, in order to prove the theorem, it is sufficient to show that the mapping $W[f, X, R]$ corresponds to mapping $w[f, X, R]$. We prove this statement by induction over complexity of $\Omega_2$-word.

If $m \in X$, $m' = R \circ m$, then, according to equations (3.17), (3.18), mappings $W[f, X, R]$ and $w[f, X, R]$ are coordinated.

Let for $m_1, ..., m_n \in X$ mappings $W[f, X, R]$ and $w[f, X, R]$ be coordinated. Let $\omega \in \Omega_2(n)$. According to the theorem 3.9

$$(3.19) \quad W(f, X, m_1 ... m_n \omega) = W(f, X, m_1) ... W(f, X, m_n) \omega$$

Because $R$ is endomorphism of $\Omega_2$-algebra $M$, then from the equation (3.19), it follows that

$$(3.20) \quad W(f, X, R \circ (m_1 ... m_n \omega)) = W(f, X, (R \circ m_1) ... (R \circ m_n) \omega)$$

$$= W(f, X, R \circ m_1) ... W(f, X, R \circ m_n) \omega$$

From equations (3.19), (3.20) and the statement of induction, it follows that the mappings $W[f, X, R]$ and $w[f, X, R]$ are coordinated for $m = m_1 ... m_n \omega$.

Let for $m_1 \in M$ mappings $W[f, X, R]$ and $w[f, X, R]$ are coordinated. Let $a \in A$. According to the theorem 3.10

$$(3.21) \quad W(f, X, am_1) = aW(f, X, m_1)$$

Because $R$ is endomorphism of representation $f$, then, from the equation (3.21), it follows that

$$(3.22) \quad W(f, X, R \circ (am_1)) = W(f, X, a R \circ m_1) = aW(f, X, R \circ m_1)$$

From equations (3.21), (3.22) and the statement of induction, it follows that mappings $W[f, X, R]$ and $w[f, X, R]$ are coordinated for $m = am_1$. $\square$

**Corollary 3.16.** *Let $X$ be the generating set of the representation $f$. Let $R$ be the endomorphism of the representation $f$. The mapping $W[f, X, R]$ is endomorphism of representation of $\Omega_1$-algebra $A$ in $\Omega_2$-algebra $W(f, X)$.* $\square$

Hereinafter we will identify mapping $W[f, X, R]$ and the set of coordinates $W(f, X, R \circ X)$.

**Theorem 3.17.** *Let $X$ be the generating set of the representation $f$. Let $R$ be the endomorphism of the representation $f$. Let $Y \subset M$. Then*

$$(3.23) \quad W(f, X, Y) \circ W(f, X, R \circ X) = W(f, X, R \circ Y)$$

$$(3.24) \quad W[f, X, R] \star W(f, X, Y) = W(f, X, R \circ Y)$$

*Proof.* The equation (3.23) follows from the equation

$$R \circ Y = (R \circ m, m \in Y)$$

as well from equations (3.13), (3.14). The equation (3.24) is corollary of equations (3.23), (3.16). □

**Theorem 3.18.** *Let $X$ be the generating set of the representation $f$. Let $R$, $S$ be the endomorphisms of the representation $f$. Then*

$$(3.25) \qquad W(f, X, S \circ X) \circ W(f, X, R \circ X) = W(f, X, R \circ S \circ X)$$

$$(3.26) \qquad W[f, X, R] \star W[f, X, S] = W[f, X, R \circ S]$$

*Proof.* The equation (3.25) follows from the equation (3.23), if we assume $Y = S \circ X$. The equation (3.26) follows from the equation (3.25) and chain of equations

$$(W[f, X, R] \star W[f, X, S]) \star W(f, X, Y)$$
$$= W[f, X, R] \star (W[f, X, S] \star W(f, X, Y))$$
$$= (W(f, X, Y) \circ W(f, X, S \circ X)) \circ W(f, X, R \circ X)$$
$$= W(f, X, Y) \circ (W(f, X, S \circ X) \circ W(f, X, R \circ X))$$
$$= W(f, X, Y) \circ W(f, X, R \circ S \circ X)$$
$$= W(f, X, R \circ S) \star W(f, X, Y)$$

□

The concept of superposition of the coordinates is very simple and resembles a kind of Turing machine. If element $m \in M$ has form either

$$m = m_1 ... m_n \omega$$

or

$$m = a m_1$$

then we are looking for the coordinates of elements $m_i$ to substitute them in an appropriate expression. As soon as an element $m \in M$ belongs to the generating set of $\Omega_2$-algebra $M$, we choose the coordinates of the corresponding element of the second factor. Therefore, we require that the second factor in the superposition has been the set of coordinates of the image of the generating set $X$.

We can generalize the definition of the superposition of coordinates and assume that one of the factors is a set of $\Omega_2$-words. Accordingly, the definition of the superposition of coordinates has the form

$$w(f, X, Y) \circ W(f, X, R \circ X) = W(f, X, Y) \circ w(f, X, R \circ X) = w(f, X, R \circ Y)$$

The following forms of writing an image of the set $Y$ under endomorphism $R$ are equivalent.

$$(3.27) \qquad R \circ Y = W(f, X, Y) \circ (R \circ X) = W(f, X, Y) \circ (W(f, X, R \circ X) \circ X)$$

From equations (3.23), (3.27), it follows that

$$(3.28) \quad (W(f, X, Y) \circ W(f, X, R \circ X)) \circ X = W(f, X, Y) \circ (W(f, X, R \circ X) \circ X)$$

The equation (3.28) is associative law for composition and allows us to write expression

$$W(f, X, Y) \circ W(f, X, R \circ X) \circ X$$

without brackets.

Consider equation ([3.25](#)), where we see change in the order of endomorphisms in a superposition of the coordinates. This equation also follows from the chain of equations, where we can immediately see when order of endomorphisms changes

$$
\begin{aligned}
&W(f, X, m) \circ W(f, X, R \circ S \circ X) \circ X \\
(3.29)\quad &= (R \circ S) \circ m = R \circ (S \circ m) \\
&= R \circ ((W(f, X, m) \circ W(f, X, S \circ X)) \circ X) \\
&= W(f, X, m) \circ W(f, X, S \circ X) \circ W(f, X, R \circ X) \circ X
\end{aligned}
$$

From the equation ([3.29](#)), it follows that coordinates of endomorphism act over coordinates of element of $\Omega_2$-algebra $M$ from the right.

**Definition 3.19.** Let $X \subset M$ be generating set of representation

$$
f : A \to {}^* M
$$

Let the mapping

$$
H : M \to M
$$

be endomorphism of the representation $f$. Let the set $X' = H \circ X$ be the image of the set $X$ under the mapping $H$. Endomorphism $H$ of representation $f$ is called **regular on generating set** $X$, if the set $X'$ is the generating set of representation $f$. Otherwise, endomorphism $H$ of representation $f$ is called **singular on generating set** $X$, $\hspace{1cm}\square$

**Definition 3.20.** Endomorphism of representation $f$ is called **regular**, if it is regular on every generating set. $\hspace{1cm}\square$

**Theorem 3.21.** *Automorphism $R$ of representation*

$$
f : A \to {}^* M
$$

*is regular endomorphism.*

*Proof.* Let $X$ be generating set of representation $f$. Let $X' = R \circ X$.

According to theorem [3.7](#) endomorphism $R$ forms the map of $\Omega_2$-words $w[f, X, R]$.

Let $m' \in M$. Since $R$ is automorphism, then there exists $m \in M$, $R \circ m = m'$. According to definition [3.6](#), $w(f, X, m)$ is $\Omega_2$-word, representing $m$ relative to generating set $X$. According to theorem [3.7](#), $w(f, X', m')$ is $\Omega_2$-word, representing of $m'$ relative to generating set $X'$

$$
w(f, X', m') = w[f, X, R](w(f, X, m))
$$

Therefore, $X'$ is generating set of representation $f$. According to definition [3.20](#), automorphism $R$ is regular. $\hspace{1cm}\square$

## 4. Basis of representation

**Definition 4.1.** If the set $X \subset M$ is generating set of representation $f$, then any set $Y$, $X \subset Y \subset M$ also is generating set of representation $f$. If there exists minimal set $X$ generating the representation $f$, then the set $X$ is called **basis of representation** $f$. $\hspace{1cm}\square$

**Theorem 4.2.** *The generating set $X$ of representation $f$ is basis iff for any $m \in X$ the set $X \setminus \{m\}$ is not generating set of representation $f$.*

*Proof.* Let $X$ be generating set of representation $f$. Assume that for some $m \in X$ there exist $\Omega_2$-word

$$(4.1) \qquad\qquad w = w(f, X \setminus \{m\}, m)$$

Consider element $m'$ such that it has $\Omega_2$-word

$$(4.2) \qquad\qquad w' = w(f, X, m')$$

that depends on $m$. According to the definition 3.6, any occurrence of $m$ into $\Omega_2$-word $w'$ can be substituted by the $\Omega_2$-word $w$. Therefore, the $\Omega_2$-word $w'$ does not depend on $m$, and the set $X \setminus \{m\}$ is generating set of representation $f$. Therefore, $X$ is not basis of representation $f$. $\qquad\square$

*Remark* 4.3. The proof of the theorem 4.2 gives us effective method for constructing the basis of the representation $f$. Choosing an arbitrary generating set, step by step, we remove from set those elements which have coordinates relative to other elements of the set. If the generating set of the representation is infinite, then this construction may not have the last step. If the representation has finite generating set, then we need a finite number of steps to construct a basis of this representation.

As noted by Paul Cohn in [4], p. 82, 83, the representation may have inequivalent bases. For instance, the cyclic group of order six has bases $\{a\}$ and $\{a^2, a^3\}$ which we cannot map one into another by endomorphism of the representation. $\qquad\square$

*Remark* 4.4. We write a basis also in following form

$$X = (x, x \in X)$$

If basis is finite, then we also use notation

$$X = (x_i, i \in I) = (x_1, ..., x_n)$$

$\qquad\qquad\qquad\qquad\qquad\qquad\qquad\qquad\qquad\qquad\qquad\qquad\square$

*Remark* 4.5. We introduced $\Omega_2$-word of $x \in M$ relative generating set $X$ in the definition 3.6. From the theorem 4.2, it follows that if the generating set $X$ is not a basis, then a choice of $\Omega_2$-word relative generating set $X$ is ambiguous. However, even if the generating set $X$ is a basis, then a representation of $m \in M$ in form of $\Omega_2$-word is ambiguous. If $m_1, ..., m_n$ are $\Omega_2$-words, $\omega \in \Omega_2(n)$ and $a \in A$, then[5]

$$(4.3) \qquad\qquad a(m_1...m_n\omega) = (am_1)...(am_n)\omega$$

It is possible that there exist equations related to specific character of representation. For instance, if $\omega$ is operation of $\Omega_1$-algebra $A$ and operation of $\Omega_2$-algebra $M$, then we require that $\Omega_2$-words $a_1...a_n\omega x$ and $a_1 x...a_n x\omega$ describe the same element of $\Omega_2$-algebra $M$.[6]

In addition to the above equations in $\Omega_2$-algebra there may be relations of the form

$$(4.4) \qquad\qquad w_1(f, X, m) = w_2(f, X, m) \quad m \notin X$$

---

[5]For instance, let $\{\overline{e}_1, \overline{e}_2\}$ be the basis of vector space over field $k$. The equation (4.3) has the form of distributive law

$$a(b^1\overline{e}_1 + b^2\overline{e}_2) = (ab^1)\overline{e}_1 + (ab^2)\overline{e}_2$$

[6]For vector space, this requirement has the form of distributive law

$$(a + b)\overline{e}_1 = a\overline{e}_1 + b\overline{e}_1$$

The feature of the equation (4.4) is that this equation cannot be reduced.[7]

On the set of $\Omega_2$-words $w(f, X)$, above equations determine **equivalence** $\rho(f)$ **generated by representation** $f$. It is evident that for any $m \in M$ the choice of appropriate $\Omega_2$-word is unique up to equivalence relations $\rho(f)$. However, if during the construction, we obtain the equality of two $\Omega_2$-word relative to given basis, then we can say without worrying about the equivalence $\rho(f)$ that these $\Omega_2$-words are equal.

A similar remark concerns the mapping $W(f, X, m)$ defined in the remark 3.8.[8]

□

**Theorem 4.6.** *Automorphism of the representation $f$ maps a basis of the representation $f$ into basis.*

*Proof.* Let the mapping $R$ be automorphism of the representation $f$. Let the set $X$ be a basis of the representation $f$. Let $X' = R \circ X$.

Assume that the set $X'$ is not basis. According to the theorem 4.2 there exists such $m' \in X'$ that $X' \setminus \{x'\}$ is generating set of the representation $f$. According to the theorem [3]-2.3.5, the mapping $R^{-1}$ is automorphism of the representation $f$. According to the theorem 3.21 and definition 3.20, the set $X \setminus \{m\}$ is generating set of the representation $f$. The contradiction completes the proof of the theorem. □

**Theorem 4.7.** *Let $X$ be the basis of the representation $f$. Let*

$$R_1 : X \to X'$$

*be arbitrary mapping of the set $X$. Consider the mapping of $\Omega_2$-words*

$$w[f, X, R_1] : w(f, X) \to w(f, X')$$

*that satisfies conditions (1), (2), (3) of the theorem 3.7 and such that*

$$x \in X => w[f, X, R_1](x) = R_1(x)$$

*There exists unique endomorphism of representation $f$[9]*

$$R : M \to M$$

*defined by rule*

$$R \circ m = w[f, X, R_1](w(f, X, m))$$

---

[7]See for instance sections [3]-5.5.2, [3]-5.5.3.

[8]If vector space has finite basis, then we represent the basis as matrix

$$\overline{\overline{e}} = \begin{pmatrix} \overline{e}_1 & ... & \overline{e}_2 \end{pmatrix}$$

We present the mapping $W(f, \overline{\overline{e}}, \overline{v})$ as matrix

$$W(f, \overline{\overline{e}}, \overline{v}) = \begin{pmatrix} v^1 \\ ... \\ v^n \end{pmatrix}$$

Then

$$W(f, \overline{\overline{e}}, \overline{v})(\overline{\overline{e}}') = W(f, \overline{\overline{e}}, \overline{v}) \begin{pmatrix} \overline{e}_1' & ... & \overline{e}_n' \end{pmatrix} = \begin{pmatrix} v^1 \\ ... \\ v^n \end{pmatrix} \begin{pmatrix} \overline{e}_1' & ... & \overline{e}_n' \end{pmatrix}$$

has form of matrix product.

[9]This statement is similar to the theorem [1]-4.1, p. 135.

*Proof.* The statement of theorem is corollary theorems 3.7, 3.11.                    □

**Corollary 4.8.** *Let $X$, $X'$ be the bases of the representation $f$. Let $R$ be the automorphism of the representation $f$ such that $X' = R \circ X$. Automorphism $R$ is uniquely defined.*                    □

*Remark* 4.9. The theorem 4.7, as well as the theorem 3.11, is the theorem of extension of mapping. However in this theorem, $X$ is not arbitrary generating set of the representation, but basis. According to remark 4.3, we cannot determine the coordinates of any element of basis through the remaining elements of the same basis. Therefore, we do not need to coordinate the mapping of the basis with representation.                    □

**Theorem 4.10.** *The set of coordinates $W(f, X, X)$ corresponds to identity transformation*

$$W[f, X, E] = W(f, X, X)$$

*Proof.* The statement of the theorem follows from the equation

$$m = W(f, X, m) \circ X = W(f, X, m) \circ W(f, X, X) \circ X$$

                    □

**Theorem 4.11.** *Let $W(f, X, R \circ X)$ be the set of coordinates of automorphism $R$. There exists set of coordinates $W(f, R \circ X, X)$, corresponding to automorphism $R^{-1}$. The set of coordinates $W(f, R \circ X, X)$ satisfy to equation*

$$(4.5) \qquad W(f, X, R \circ X) \circ W(f, R \circ X, X) = W(f, X, X)$$

$$W[f, X, R^{-1}] = W[f, X, R]^{-1} = W(f, R \circ X, X)$$

*Proof.* Since $R$ is automorphism of the representation $f$, then, according to the theorem 4.6, the set $R \circ X$ is a basis of the representation $f$. Therefore, there exists the set of coordinates $W(f, R \circ X, X)$. The equation (4.5) follows from the chain of equations

$$W(f, X, R \circ X) \circ W(f, R \circ X, X) = W(f, X, R \circ X) \circ W(f, X, R^{-1} \circ X)$$
$$= W(f, X, R^{-1} \circ R \circ X) = W(f, X, X)$$

                    □

*Remark* 4.12. In $\Omega_2$-algebra $M$ there is no universal algorithm for determining the set of coordinates $W(f, R \circ X, X)$ for given set $W(f, X, R \circ X)$.[10] We assume that in the theorem 4.11 this algorithm is given implicitly. It is evident also that the set of $\Omega_2$-words

$$(4.6) \qquad W(f, X, R \circ X) \circ W(f, R \circ X, X) \circ X$$

in general, does not coincide with the set of $\Omega_2$-words

$$(4.7) \qquad W(f, X, X) \circ X$$

The theorem 4.11 states that sets of $\Omega_2$-words (4.6) and (4.7) coincide up to equivalence generated by the representation $f$.                    □

---

[10]In vector space, the matrix of numbers corresponds to linear transformation. Accordingly, the inverse matrix corresponds to inverse transformation.

**Theorem 4.13.** *The group of automorphisms $G(f)$ of effective representation $f$ in $\Omega_2$-algebra $M$ generates effective representation in $\Omega_2$-algebra $M$.*

*Proof.* From the corollary 4.8, it follows that if automorphism $R$ maps a basis $X$ into a basis $X'$, then the set of coordinates $W(f, X, X')$ uniquely determines an automorphism $R$. From the theorem 3.15, it follows that the set of coordinates $W(f, X, X')$ determines the mapping of coordinates relative to the basis $X$ under automorphism of the representation $f$. From the equation (3.27), it follows that automorphism $R$ acts from the right on elements of $\Omega_2$-algebra $M$. From the equation (3.25), it follows that the representation of group covariant from right. According to the theorem 4.10, the set of coordinates $W(f, X, X)$ corresponds to identity transformation. From the theorem 4.11, it follows that the set of coordinates $W(f, R \circ X, X)$ corresponds to transformation, inverse to transformation $W(f, X, R \circ X)$. $\qquad\square$

## 5. Basis Manifold of Representation

The set $\mathcal{B}(f)$ of bases of representation $f$ is called **basis manifold of representation** $f$.

According to theorem 4.6, automorphism $R$ of the representation $f$ generates transformation

(5.1)
$$R : Y \to R \circ Y$$
$$W(f, X, R \circ Y) = W(f, X, Y) \circ W(f, X, R \circ X)$$

of the basis manifold of representation. This transformation is called **active**. According to the theorem [3]-2.3.5, we defined **active representation**

$$A(f) : G(f) \to \mathcal{B}^*(f)$$

of group $G(f)$ in basis manifold $\mathcal{B}(f)$. According to the corollary 4.8, this representation is single transitive.

*Remark* 5.1. According to remark 4.3, it is possible that there exist bases such that there is no active transformation between them. Then we consider the orbit of selected basis as basis manifold. Therefore, it is possible that the representation $f$ has different basis manifolds. We will assume that we have chosen a basis manifold. $\qquad\square$

**Definition 5.2.** Automorphism $S$ of representation $A(f)$ is called **passive transformation of the basis manifold of representation**. We also use notation

$$S(X) = S \star X$$

to denote the image of basis $X$ under passive transformation $S$. $\qquad\square$

**Theorem 5.3.** *Passive transformation of basis has form*

(5.2)
$$S : Y \to S \star Y$$
$$W(f, X, S \star Y) = W(f, X, S \star Y) \circ W(f, X, Y)$$

*Proof.* According to the equation (5.1), active transformation acts from right on coordinates of basis. The equation (5.2) follows from theorems [3]-2.5.12, [3]-2.5.13, [3]-2.5.15, according to these theorems, passive transformation acts from left on coordinates of basis. $\qquad\square$

**Theorem 5.4.** *Let $S$ be passive transformation of the basis manifold of the representation $f$. Let $X$ be the basis of the representation $f$, $X' = S \star X$. For basis $Y$, let there exists an active transformation $R$ such that $Y = R \circ X$. Assume $Y' = R \circ X'$. Then $S \star Y = Y'$.*

*Proof.* According to equation (5.1), active transformation of coordinates of basis $Y$ has form

$$R : X' = S \star X \to Y'$$
(5.3)
$$W(f, X, Y') = W(f, X, X') \circ W(f, X, Y)$$

Let $Y'' = S \star Y$. From the equation (5.2), it follows that

$$S : Y = R \circ X \to Y''$$
(5.4)
$$W(f, X, Y'') = W(f, X, X') \circ W(f, X, Y)$$

From match of expressions in equations (5.3), (5.4), it follows that $Y' = Y'' = S \star Y$. Therefore, the diagram

is commutative. □

**Theorem 5.5.** *There exists single transitive* **passive representation**

$$P(f) : G(f) \to {}^*\mathcal{B}(f)$$

*Proof.* Since $A(f)$ is single transitive representation of group $G(f)$, then according to the definition 5.2 and the theorem [3]-2.3.5, the set of passive transformations forms group. □

Let passive transformation $S$ maps basis $Y$ into basis

$$Y' = S \star Y$$
(5.5)
$$W(f, X, Y') = W(f, X, S \star X) \circ W(f, X, Y)$$

Let passive transformation $T$ maps basis $Y'$ into basis

$$Y'' = T \star Y'$$
(5.6)
$$W(f, X, Y'') = W(f, X, T \star X) \circ W(f, X, Y')$$

Denote $T \star S$ the product of passive transformations $T$ and $S$

$$Y'' = T \star S \star Y$$
(5.7)
$$W(f, X, Y'') = W(f, X, T \star S \star X) \circ W(f, X, Y)$$

From equations (5.5), (5.6), it follows that

$$Y'' = T \star S \star Y$$
(5.8)
$$W(f, X, Y'') = W(f, X, T \star X) \circ W(f, X, S \star X) \circ W(f, X, Y)$$

From equations (5.7), (5.8) and theorem 5.5, it follows that

$$(5.9) \qquad W(f, X, T \star S \star X) = W(f, X, T \star X) \circ W(f, X, S \star X)$$

## 6. Geometric Object of Representation of Universal Algebra

An active transformation changes bases and elements of $\Omega_2$-algebra uniformly and coordinates of element relative basis do not change. A passive transformation changes only the basis and it leads to change of coordinates of element relative to the basis.

Let passive transformation $S \in G(f)$ maps basis $Y \in \mathcal{B}(f)$ into basis $Y' \in \mathcal{B}(f)$

$$(6.1) \qquad W(f, X, Y') = W(f, X, S \star X) \circ W(f, X, Y)$$

Let element $m \in M$ has $\Omega_2$-word

$$(6.2) \qquad m = W(f, Y, m) \circ Y$$

relative to basis $Y$ and has $\Omega_2$-word

$$(6.3) \qquad m = W(f, Y', m) \circ Y'$$

relative to basis $Y'$. From (6.1) and (6.3), it follows that

$$
\begin{aligned}
(6.4) \qquad W(f, Y, m) \circ Y &= W(f, Y', m) \circ Y' = W(f, Y', m) \circ W(f, Y, Y') \circ Y \\
&= W(f, Y', m) \circ W(f, Y, S \star Y) \circ Y
\end{aligned}
$$

Comparing (6.2) and (6.4) we get

$$(6.5) \qquad W(f, Y, m) = W(f, Y', m) \circ W(f, Y, S \star Y)$$

Because $S$ is automorphism we get from (6.5) and the theorem 4.11

$$(6.6) \qquad W(f, Y', m) = W(f, Y, m) \circ W(f, S \star Y, Y)$$

Coordinate transformation (6.6) does not depend on element $m$ or basis $Y$, but is defined only by coordinates of element $m$ relative to basis $Y$.

**Theorem 6.1.** *Coordinate transformations* (6.6) *form contravariant effective $T\star$-representation of group $G(f)$ which is called* **coordinate representation in $\Omega_2$-algebra**.

*Proof.* According to corollary 3.16, the transformation (6.6) is the endomorphism of representation of $\Omega_1$-algebra $A$ into $\Omega_2$-algebra $W(f, X)$.[11]

Suppose we have two consecutive passive transformations $S$ and $T$. Coordinate transformation (6.6) corresponds to passive transformation $S$. Coordinate transformation

$$(6.7) \qquad W(f, Y'', m) = W(f, Y', m) \circ W(f, T \star Y, Y)$$

corresponds to passive transformation $T$. Product of coordinate transformations (6.6) and (6.7) has form

$$
\begin{aligned}
(6.8) \qquad W(f, Y'', m) &= W(f, Y, m) \circ W(f, S \star Y, Y) \circ W(f, T \star Y, Y) \\
&= W(f, Y, m) \circ W(f, T \star S \star Y, Y)
\end{aligned}
$$

---

[11] This transformation does not generate an endomorphism of the representation $f$. Coordinates change because basis relative which we determinate coordinates changes. However, element of $\Omega_2$-algebra $M$, coordinates of which we are considering, does not change.

and is coordinate transformation corresponding to passive transformation $T \star S$. It proves that coordinate transformations form contravariant $T\star$-representation of group $G(f)$.

Suppose coordinate transformation does not change coordinates of selected basis. Then unit of group $G(f)$ corresponds to it because representation is single transitive. Therefore, coordinate representation is effective. $\qquad\square$

Let $f$ be representation of $\Omega_1$-algebra $A$ in $\Omega_2$-algebra $M$. Let $g$ be representation of $\Omega_1$-algebra $A$ in $\Omega_3$-algebra $N$. Passive representation $g_p$ coordinated with passive representation $f_p$, if there exists homomorphism $h$ of group $G(f)$ into group $G(g)$. Consider diagram[12]

$$
\begin{array}{ccc}
{}^*\mathcal{B}(f) & \xrightarrow{\;H\;} & {}^*\mathcal{B}(g) \\[2pt]
P(f)\Big\uparrow & \nearrow{\scriptstyle f'} & \Big\uparrow P(g) \\[2pt]
G(f) & \xrightarrow{\;h\;} & G(g)
\end{array}
$$

Since mappings $P(f)$, $P(g)$ are isomorphisms of group, then mapping $H$ is homomorphism of groups. Therefore, mapping $f'$ is representation of group $G(f)$ in basis manifold $\mathcal{B}(g)$. According to design, passive transformation $H(S)$ of $\Omega_3$-algebra $N$ corresponds to passive transformation $S$ of $\Omega_2$-algebra $M$.

(6.9)                         $Y'_N = H(S) \star Y_N$

Then coordinate transformation in $N$ gets form

(6.10)          $W(g, Y'_P, m) = W(g, Y_N, m) \circ W(g, H(S) \star Y_P, Y_P)$

**Definition 6.2.** Orbit

$$(W(g, Y_N, n) \circ W(g, H(S) \star Y_N, Y_N), S \star Y_M, S \in G(f))$$

is called **geometric object in coordinate representation defined in $\Omega_2$-algebra** $M$. For any basis $Y'_M = S \star Y_M$ corresponding point (6.10) of orbit defines **coordinates of geometric object in coordinate space of representation** relative basis $Y'_M$. $\qquad\square$

**Definition 6.3.** Orbit

$$(W(g, Y_N, m) \circ W(g, H(S) \star Y_N, Y_N), H(S) \star Y_N, S \star Y_M, S \in G(f))$$

is called **geometric object defined in $\Omega_2$-algebra** $M$. For any basis $Y'_M = S \star Y_M$ corresponding point (6.10) of orbit defines **coordinates of a geometric object in $\Omega_2$-algebra** $M$ relative to basis $Y'_M$ and the corresponding vector

$$p = W(g, Y'_P, p) \circ Y'_P$$

is called **representative of geometric object in $\Omega_2$-algebra** $M$. $\qquad\square$

Since a geometric object is an orbit of representation, we see that according to the theorem [3]-2.4.13 the definition of the geometric object is a proper definition.

We also say that $p$ is a **geometric object of type** $H$

---

[12]We can relax definition of consistency of representations and assume that $N$ is $\Omega_2$-algebra. Then it is sufficient to require that the mapping $(h, H)$ is morphism of representations. However in this case, elements of affine representation cannot be geometric object of vector space.

Definition [6.2] introduces a geometric object in coordinate space. We assume in definition [6.3] that we selected a basis of representation $g$. This allows using a representative of the geometric object instead of its coordinates.

**Theorem 6.4** (**invariance principle**). *Representative of geometric object does not depend on selection of basis $Y_M'$.*

*Proof.* To define representative of geometric object, we need to select basis $Y_M$, basis $Y_P$ and coordinates of geometric object $W(g, Y_P, p)$. Corresponding representative of geometric object has form

$$p = W(g, Y_P, p) \circ Y_P$$

Suppose we map basis $Y_M$ to basis $Y_M'$ by passive transformation

$$Y_M' = S \circ Y_M$$

According building this forms passive transformation [(6.9)] and coordinate transformation [(6.10)]. Corresponding representative of geometric object has form

$$\begin{aligned} p' &= W(g, Y_P', p') \circ Y_P' \\ &= W(g, Y_P, p) \circ W(g, H(S) \star Y_P, Y_P) \circ W(g, Y_P, H(S) \star Y_P) \circ Y_P \\ &= W(g, Y_P, p) \circ Y_P = p \end{aligned}$$

Therefore representative of geometric object is invariant relative selection of basis. $\square$

**Theorem 6.5.** *The set of geometric objects of type $H$ is $\Omega_3$-algebra.*

*Proof.* Let

$$p_i = W(g, Y_P, p_i) \circ Y_P \quad i = 1, ..., n$$

For operation $\omega \in \Omega_3(n)$ we assume

(6.11) $$p_1 ... p_n \omega = W(g, Y_P, p_1) ... W(g, Y_P, p_n) \omega \circ Y_P$$

Since for arbitrary endomorphism $S$ of $\Omega_2$-algebra $M$, the mapping $W(g, H(S) \star Y_N, Y_N)$ is endomorphism of $\Omega_3$-algebra $N$, then the definition [(6.11)] is correct. $\square$

**Theorem 6.6.** *There exists the representation of $\Omega_1$-algebra $A$ in $\Omega_3$-algebra $N$ of geometric objects of type $H$.*

*Proof.* Let

$$p = W(g, Y_P, p) \circ Y_P$$

For $a \in A$, we assume

(6.12) $$ap = aW(g, Y_P, p_1) \circ Y_P$$

Since for arbitrary endomorphism $S$ of $\Omega_2$-algebra $M$, the mapping $W(g, H(S) \star Y_N, Y_N)$ is endomorphism of representation $g$, then the definition [(6.12)] is correct. $\square$

## 7. Examples of Basis of Representation of Universal Algebra

7.1. **Vector Space.** Consider the vector space $\overline{V}$ over the field $F$. Given the set of vectors $\overline{e}_1$, ..., $\overline{e}_n$, according to algorithm of construction of coordinates over vector space, coordinates include such elements as $\overline{e}_1 + \overline{e}_2$ and $a\overline{e}_1$. Recursively using rules, contained in the definition 3.6, we conclude that the set of vectors $\overline{e}_1$, ..., $\overline{e}_n$, generates the set of linear combinations

$$a^1\overline{e}_1 + ... + a^n\overline{e}_n$$

According to the theorem 4.2, the set of vectors $\overline{e}_1$, ..., $\overline{e}_n$, is a basis if for any $i$, $i = 1$, ..., $n$, vector $\overline{e}_i$ is not linear combination of other vectors. This requirement is equivalent to the requirement of linear independence of vectors.

7.2. **Representation of Group on the Set.** Consider the representation from the example [3]-2.6.2. We can consider the set $M$ as union of orbits of the representation of the group $G$. We can select for basis of the representation the set of points such that one and only one point belongs to each orbit. If $X$ is the basis of representation, $A \in X$, $g \in G$, then $\Omega_2$-word has form $A + g$. Since there is no operations on the set $M$, then there is no $\Omega_2$-word containing different elements of the basis. If representation of group $G$ is single transitive, then basis of representation consists of one point. Any point of the set $M$ can be such point.

## 9. Index

## 10. Special Symbols and Notations

$A(f)$   active representation of group $G(f)$
in basis manifold $\mathcal{B}(f)$ **15**

$\mathcal{B}(f)$   basis manifold of representation $f$ **15**

$\mathcal{B}_f$   lattice of subrepresentations of
representation $f$ **3**

$J(f)$   closure operator of representation $f$ **3**

$P(f)$   passive representation of group $G(f)$
in basis manifold $\mathcal{B}(f)$ **16**

$R \circ m$   image of $m$ under endomorphism $R$
of effective representation **4**

$S \star X$   image of basis $X$ under passive
transformation $S$ **15**

$W(f, X)$   set of coordinates of
representation $J(f, X)$ **6**

$W(f, X, m)$   coordinates of element $m$ of
representation $f$ relative to set $X$ **5**

$(W(g, Y_N, n) \circ W(g, H(S) \star Y_N, Y_N), S \star Y_M, S \in G(f))$
geometric object in coordinate
representation defined in $\Omega_2$-algebra
$M$ **18**

$(W(g, Y_N, m) \circ W(g, H(S) \star Y_N, Y_N), H(S) \star Y_N, S \star Y_M, S \in G(f))$
geometric object defined in $\Omega_2$-algebra
$M$ **18**

$W(f, X, B)$   set of coordinates of set
$B \subset J(f, X)$ **6**

$w(f, X, B)$   set of $\Omega_2$-words representing
set $B \subset J(f, X)$ **4**

$W(f, X, m) \circ W(f, X, R \circ X)$   superposition
of coordinates of the representation $f$
and the element $m$ **8**

$w(f, X, m)$   $\Omega_2$-word representing element
$m \in J(f, X)$ **4**

$w(f, X)$   set of $\Omega_2$-words of representation
$J(f, X)$ **5**

# Базис представления универсальной алгебры

Александр Клейн

**Аннотация.** Если множество эндоморфизмов универсальной алгебры $A$ имеет структуру универсальной алгебры $B$, то мы говорим, что определено представление универсальной алгебры $B$ в универсальной алгебре $A$. Следовательно, роль представления универсальной алгебры аналогична роли симметрии в геометрии и физике. Морфизм представления - это отображение, сохраняющее структуру представления. Изучение морфизмов представлений ведёт к понятиям множества образующих и базиса представления. Множество автоморфизмов представления универсальной алгебры порождает группу. Парные представления этой группы в многообразии базисов представления называются активным и пассивным представлениями. Пассивное представление является основой понятия геометрического объекта и теории инвариантов представления универсальной алгебры.

## Содержание

## 1. Предисловие

Эта статья написана на основе главы [3]-3.

Роль представления универсальной алгебры аналогична роли симметрии в геометрии и физике В обоих случаях мы изучаем структурированные множества преобразований, и знание структуры этого множества даёт нам дополнительные сведения об изучаемом объекте.

Aleks_Kleyn@MailAPS.org.
http://sites.google.com/site/AleksKleyn/.
http://arxiv.org/a/kleyn_a_1.
http://AleksKleyn.blogspot.com/.

Мы рассматриваем теорию представлений универсальной алгебры как расширение теории универсальной алгебры. Любая алгебраическая структура предполагает множество отображений, сохраняющих эту структуру. Отображение, сохраняющее структуру представления универсальной алгебры, называется морфизмом представлений.

Частным случаем морфизма представлений является автоморфизм представления. Изучение автоморфизмов представления непосредственно связано с необходимостью ответить на вопрос какова структура множества образующих представления. Если представление имеет минимальное множество образующих, то такое множество называется базисом представления. Автоморфизм представления отображает базис в базис.

Множество автоморфизмов представления порождает группу. Парные представления этой группы в многообразии базисов представления называются активным и пассивным представлениями. Пассивное представление является основой понятия геометрического объекта и теории инвариантов представления универсальной алгебры.

## 2. Соглашения

**Соглашение 2.1.** *В [4] произвольная операция алгебры обозначена буквой $\omega$, и $\Omega$ - множество операций некоторой универсальной алгебры. Соответственно, универсальная алгебра с множеством операций $\Omega$ обозначается $\Omega$-алгебра. Аналогичные обозначения мы видим в [2] с той небольшой разницей, что операция в алгебре обозначена буквой $f$ и $\mathcal{F}$ - множество операций. Я выбрал первый вариант обозначений, так как в этом случае легче видно, где я использую операцию.* □

**Соглашение 2.2.** *Пусть $A$ - $\Omega_1$-алгебра. Пусть $B$ - $\Omega_2$-алгебра. Запись*

$$A \longrightarrow_* B$$

*означает, что определено представление $\Omega_1$-алгебры $A$ в $\Omega_2$-алгебре $B$.* □

Без сомнения, у читателя могут быть вопросы, замечания, возражения. Я буду признателен любому отзыву.

## 3. Множество образующих представления

**Определение 3.1.** Пусть

$$f : A \rightarrow {}^*M$$

представление $\Omega_1$-алгебры $A$ в $\Omega_2$-алгебре $M$. Множество $N \subset M$ называется **стабильным множеством представления** $f$, если $f(a)(m) \in N$ для любых $a \in A$, $m \in N$. □

Мы также будем говорить, что множество $M$ стабильно относительно представления $f$.

**Теорема 3.2.** *Пусть*

$$f : A \rightarrow {}^*M$$

*представление $\Omega_1$-алгебры $A$ в $\Omega_2$-алгебре $M$. Пусть множество $N \subset M$ является подалгеброй $\Omega_2$-алгебры $M$ и стабильным множеством представления $f$. Тогда существует представление*

$$f_N : A \rightarrow {}^*N$$

*такое, что* $f_N(a) = f(a)|_N$. *Представление* $f_N$ *называется* **подпредставлением представления** $f$.

*Доказательство.* Пусть $\omega_1$ - $n$-арная операция $\Omega_1$-алгебры $A$. Тогда для любых $a_1, ..., a_n \in A$ и любого $b \in N$

$$(f_N(a_1)...f_N(a_n)\omega_1)(b) = (f(a_1)...f(a_n)\omega_1)(b)$$
$$= f(a_1...a_n\omega_1)(b)$$
$$= f_N(a_1...a_n\omega_1)(b)$$

Пусть $\omega_2$ - $n$-арная операция $\Omega_2$-алгебры $M$. Тогда для любых $b_1, ..., b_n \in N$ и любого $a \in A$

$$f_N(a)(b_1)...f_N(a)(b_n)\omega_2 = f(a)(b_1)...f(a)(b_n)\omega_2$$
$$= f(a)(b_1...b_n\omega_2)$$
$$= f_N(a)(b_1...b_n\omega_2)$$

Утверждение теоремы доказано. $\square$

Из теоремы 3.2 следует, что если $f_N$ - подпредставление представления $f$, то отображение $(id : A \to A, id_N : N \to M)$ является морфизмом представлений.

**Теорема 3.3.** *Множество[1] $\mathcal{B}_f$ всех подпредставлений представления $f$ порождает систему замыканий на $\Omega_2$-алгебре $M$ и, следовательно, является полной структурой.*

*Доказательство.* Пусть $(K_\lambda)_{\lambda \in \Lambda}$ - семейство подалгебр $\Omega_2$-алгебры $M$, стабильных относительно представления $f$. Операцию пересечения на множестве $\mathcal{B}_f$ мы определим согласно правилу

$$\bigcap f_{K_\lambda} = f_{\cap K_\lambda}$$

Операция пересечения подпредставлений определена корректно. $\cap K_\lambda$ - подалгебра $\Omega_2$-алгебры $M$. Пусть $m \in \cap K_\lambda$. Для любого $\lambda \in \Lambda$ и для любого $a \in A$, $f(a)(m) \in K_\lambda$. Следовательно, $f(a)(m) \in \cap K_\lambda$. Следовательно, $\cap K_\lambda$ - стабильное множество представления $f$. $\square$

Обозначим соответствующий оператор замыкания через $J(f)$. Таким образом, $J(f, X)$ является пересечением всех подалгебр $\Omega_2$-алгебры $M$, содержащих $X$ и стабильных относительно представления $f$.

**Теорема 3.4.** *Пусть[2]*

$$f : A \to {}^*M$$

*представление $\Omega_1$-алгебры $A$ в $\Omega_2$-алгебре $M$. Пусть $X \subset M$. Определим подмножество $X_k \subset M$ индукцией по $k$.*

$$X_0 = X$$
$$x \in X_k => x \in X_{k+1}$$
$$x_1 \in X_k, ..., x_n \in X_k, \omega \in \Omega_2(n) => x_1...x_n\omega \in X_{k+1}$$
$$x \in X_k, a \in A => f(a)(x) \in X_{k+1}$$

---

[1]Это определение аналогично определению структуры подалгебр ([4], стр. 93, 94)

[2]Утверждение теоремы аналогично утверждению теоремы 5.1, [4], стр. 94.

*Тогда*

$$\bigcup_{k=0}^{\infty} X_k = J(f, X)$$

*Доказательство.* Если положим $U = \cup X_k$, то по определению $X_k$ имеем $X_0 \subset J(f, X)$, и если $X_k \subset J(f, X)$, то $X_{k+1} \subset J(f, X)$. По индукции следует, что $X_k \subset J(f, X)$ для всех $k$. Следовательно,

$$(3.1) \qquad\qquad U \subset J(f, X)$$

Если $a \in U^n$, $a = (a_1, ..., a_n)$, где $a_i \in X_{k_i}$, и если $k = \max\{k_1, ..., k_n\}$, то $a_1...a_n\omega \in X_{k+1} \subset U$. Следовательно, $U$ является подалгеброй $\Omega_2$-алгебры $M$.

Если $m \in U$, то $m \in X_k$ для некоторого $k$. Следовательно, $f(a)(m) \in X_{k+1} \subset U$ для любого $a \in A$. Следовательно, $U$ - стабильное множество представления $f$.

Так как $U$ - подалгеброй $\Omega_2$-алгебры $M$ и стабильное множество представления $f$, то определено подпредставление $f_U$. Следовательно,

$$(3.2) \qquad\qquad J(f, X) \subset U$$

Из (3.1), (3.2), следует $J(f, X) = U$.                        $\square$

**Определение 3.5.** $J(f, X)$ называется **подпредставлением, порождённым множеством** $X$, а $X$ - **множеством образующих подпредставления** $J(f, X)$. В частности, **множеством образующих представления** $f$ будет такое подмножество $X \subset M$, что $J(f, X) = M$.

Нетрудно видеть, что определение множества образующих представления не зависит от того, эффективно представление или нет. Поэтому в дальнейшем мы будем предполагать, что представление эффективно и будем опираться на соглашение для эффективного $T\star$-представления в замечании [3]-2.1.9. Мы также будем пользоваться записью

$$R \circ m = R(m)$$

для образа $m \in M$ при эндоморфизме эффективного представления. Согласно определению произведения отображений, для любых эндоморфизмов $R$, $S$ верно равенство

$$(3.3) \qquad\qquad (R \circ S) \circ m = R \circ (S \circ m)$$

Равенство (3.3) является законом ассоциативности для $\circ$ и позволяет записать выражение

$$R \circ S \circ m$$

без использования скобок.

Из теоремы 3.4 следует следующее определение.

**Определение 3.6.** Пусть $X \subset M$. Для любого $x \in J(f, X)$ существует $\Omega_2$-слово, определённое согласно следующему правилам.

(1) Если $m \in X$, то $m$ - $\Omega_2$-слово.
(2) Если $m_1, ..., m_n$ - $\Omega_2$-слова и $\omega \in \Omega_2(n)$, то $m_1...m_n\omega$ - $\Omega_2$-слово.
(3) Если $m$ - $\Omega_2$-слово и $a \in A$, то $am$ - $\Omega_2$-слово.

$\Omega_2$-**слово** $w(f, X, m)$ представляет данный элемент $m \in J(f, X)$. Мы будем отождествлять элемент $m \in J(f, X)$ и соответствующее ему $\Omega_2$-слово, выражая это равенством

$$m = w(f, X, m)$$

Аналогично, для произвольного множества $B \subset J(f, X)$ рассмотрим множество $\Omega_2$-слов[3]

$$w(f, X, B) = \{w(f, X, m) : m \in B\}$$

Мы будем также пользоваться записью

$$w(f, X, B) = (w(f, X, m), m \in B)$$

Обозначим $w(f, X)$ **множество $\Omega_2$-слов представления** $J(f, X)$. □

**Теорема 3.7.** *Эндоморфизм $R$ представления*

$$f : A \to {}^*M$$

*порождает отображение $\Omega_2$-слов*

$$w[f, X, R] : w(f, X) \to w(f, X') \quad X \subset M \quad X' = R \circ X$$

*такое, что*

(1) *Если $m \in X$, $m' = R \circ m$, то*

$$w[f, X, R](m) = m'$$

(2) *Если*

$$m_1, ..., m_n \in w(f, X)$$

$$m'_1 = w[f, X, R](m_1) \quad ... \quad m'_n = w[f, X, R](m_n)$$

*то для операции $\omega \in \Omega_2(n)$ справедливо*

$$w[f, X, R](m_1...m_n\omega) = m'_1...m'_n\omega$$

(3) *Если*

$$m \in w(f, X) \quad m' = w[f, X, R](m) \quad a \in A$$

*то*

$$w[f, X, R](am) = am'$$

*Доказательство.* Утверждения (1), (2) теоремы справедливы в силу определения эндоморфизма $R$. Утверждение (3) теоремы следует из равенства [3]-(2.2.4), так как $r = \mathrm{id}$. □

---

[3]Выражение $w(f, X, m)$ является частным случаем выражения $w(f, X, B)$, а именно

$$w(f, X, \{m\}) = \{w(f, X, m)\}$$

*Замечание* 3.8. Пусть $R$ - эндоморфизм представления $f$. Пусть

$$m \in J(f, X) \quad m' = R \circ m \quad X' = R \circ X$$

Теорема 3.7 утверждает, что $m' \in J_f(X')$ . Теорема 3.7 также утверждает, что $\Omega_2$-слово, представляющее $m$, относительно $X$ и $\Omega_2$-слово, представляющее $m'$, относительно $X'$ формируются согласно одному и тому же алгоритму. Это позволяет рассматривать множество $\Omega_2$-слов $w(f, X', m')$ как отображение

$$W(f, X, m) : X' \to w(f, X', m')$$

$$W(f, X, m)(X') = W(f, X, m) \circ X'$$

такое, что, если для некоторого эндоморфизма $R$

$$X' = R \circ X \quad m' = R \circ m$$

то

$$W(f, X, m) \circ X' = w(f, X', m') = m'$$

Отображение $W(f, X, m)$ называется **координатами элемента** $m$ **относительно множества** $X$. Аналогично, мы можем рассмотреть координаты множества $B \subset J(f, X)$ относительно множества $X$

$$W(f, X, B) = \{W(f, X, m) : m \in B\} = (W(f, X, m), m \in B)$$

Обозначим $W(f, X)$ **множество координат представления** $J(f, X)$. $\quad\square$

**Теорема 3.9.** *На множестве координат $W(f, X)$ определена структура $\Omega_2$-алгебры.*

*Доказательство.* Пусть $\omega \in \Omega_2(n)$. Тогда для любых $m_1, ..., m_n \in J(f, X)$ положим

$$(3.4) \qquad W(f, X, m_1)...W(f, X, m_n)\omega = W(f, X, m_1...m_n\omega)$$

Согласно замечанию 3.8, из равенства (3.4) следует

$$(3.5) \qquad (W(f, X, m_1)...W(f, X, m_n)\omega) \circ X = W(f, X, m_1...m_n\omega) \circ X$$
$$= w(f, X, m_1...m_n\omega)$$

Согласно правилу (2) определения 3.6, из равенства (3.5) следует
$$(3.6)$$
$$(W(f, X, m_1)...W(f, X, m_n)\omega) \circ X = w(f, X, m_1)...w(f, X, m_n)\omega$$
$$= (W(f, X, m_1) \circ X)...(W(f, X, m_n) \circ X)\omega$$

Из равенства (3.6) следует корректность определения (3.4) операции $\omega$ на множестве координат. $\quad\square$

**Теорема 3.10.** *Определено представление $\Omega_1$-алгебры $A$ в $\Omega_2$-алгебре $W(f, X)$.*

*Доказательство.* Пусть $a \in A$. Тогда для любого $m \in J(f, X)$ положим
$$(3.7) \qquad aW(f, X, m) = W(f, X, am)$$

Согласно замечанию 3.8, из равенства (3.7) следует
$$(3.8) \qquad (aW(f, X, m)) \circ X = W(f, X, am) \circ X = w(f, X, am)$$

Согласно правилу (3) определения 3.6, из равенства (3.8) следует

$$(3.9) \qquad (aW(f, X, m)) \circ X = aw(f, X, m) = a(W(f, X, m) \circ X)$$

Из равенства (3.9) следует корректность определения (3.7) представления $\Omega_1$-алгебры $A$ в $\Omega_2$-алгебре $W(f, X)$.          $\square$

**Теорема 3.11.** *Пусть*

$$f : A \to {}^*M$$

*представление $\Omega_1$-алгебры $A$ в $\Omega_2$-алгебре $M$. Для заданных множеств $X \subset M$, $X' \subset M$, пусть отображение*

$$R_1 : X \to X'$$

*согласовано со структурой представления $f$, т. е.*

$$\omega \in \Omega_2(n) \ x_1, \ ..., \ x_n, \ x_1...x_n\omega \in X, \ R_1(x_1...x_n\omega) \in X'$$

$$=> R_1(x_1...x_n\omega) = R_1(x_1)...R_1(x_n)\omega$$

$$x \in X, \ a \in A, \ R_1(ax) \in X'$$

$$=> R_1(ax) = aR_1(x)$$

*Рассмотрим отображение $\Omega_2$-слов*

$$w[f, X, R_1] : w(f, X) \to w(f, X')$$

*удовлетворяющее условиям (1), (2), (3) теоремы 3.7, и такое, что*

$$x \in X => w[f, X, R_1](x) = R_1(x)$$

*Существует единственный эндоморфизм $\Omega_2$-алгебры $M$*

$$R : M \to M$$

*определённый правилом*

$$R \circ m = w[f, X, R_1](w(f, X, m))$$

*который является морфизмом представлений $J(f, X)$ и $J(f, X')$.*

*Доказательство.* Мы будем доказывать теорему индукцией по сложности $\Omega_2$-слова.

Если $w(f, X, m) = m$, то $m \in X$. Согласно условию (1) теоремы 3.7,

$$R \circ m = w[f, X, R_1](w(f, X, m)) = w[f, X, R_1](m) = R_1(m)$$

Следовательно, на множестве $X$ отображения $R$ и $R_1$ совпадают, и отображение $R$ согласовано со структурой представления $f$.

Пусть $\omega \in \Omega_2(n)$. Пусть отображение $R$ определено для $m_1, ..., m_n \in J(f, X)$. Пусть

$$w_1 = w(f, X, m_1) \quad ... \quad w_n = w(f, X, m_n)$$

Если $m = m_1...m_n\omega$, то согласно правилу (2) определения 3.6,

$$w(f, X, m) = w_1...w_n\omega$$

Согласно условию (2) теоремы 3.7,

$$R \circ m = w[f, X, R_1](w(f, X, m)) = w[f, X, R_1](w_1...w_n\omega)$$
$$= w[f, X, R_1](w_1)...w[f, X, R_1](w_n)\omega$$
$$= (R \circ m_1)...(R \circ m_n)\omega$$

Следовательно, отображение $R$ является эндоморфизмом $\Omega_2$-алгебры $M$.

Пусть отображение $R$ определено для $m_1 \in J(f, X)$, $w_1 = w(f, X, m_1)$. Пусть $a \in A$. Если $m = am_1$, то согласно правилу (3) определения 3.6,

$$w(f, X, am_1) = aw_1$$

Согласно условию (3) теоремы 3.7,

$$R \circ m = w[f, X, R_1](w(f, X, m)) = w[f, X, R_1](aw_1)$$
$$= aw[f, X, R_1](w_1) = aR \circ m_1$$

Из равенства [3]-(2.2.4) следует, что отображение $R$ является морфизмом представления $f$.

Единственность эндоморфизма $R$, а следовательно, корректность его определения, следует из следующего рассуждения. Допустим, $m \in M$ имеет различные $\Omega_2$-слова относительно множества $X$, например

(3.10) $$m = x_1...x_n\omega = ax$$

Так как $R$ - эндоморфизм представления, то из равенства (3.10) следует

(3.11) $$R \circ m = R \circ (x_1...x_n\omega) = (R \circ x_1)...(R \circ x_n)\omega = R \circ (ax) = aR \circ x$$

Из равенства (3.11) следует

(3.12) $$R \circ m = (R \circ x_1)...(R \circ x_n)\omega = aR \circ x$$

Из равенств (3.10), (3.12) следует, что равенство (3.10) сохраняется при отображении. Следовательно, образ $m$ не зависит от выбора координат. $\qquad\square$

*Замечание* 3.12. Теорема 3.11 - это теорема о продолжении отображения. Единственное, что нам известно о множестве $X$ - это то, что $X$ - множество образующих представления $f$. Однако, между элементами множества $X$ могут существовать соотношения, порождённые либо операциями $\Omega_2$-алгебры $M$, либо преобразованиями представления $f$. Поэтому произвольное отображение множества $X$, вообще говоря, не может быть продолжено до эндоморфизма представления $f$.[4] Однако, если отображение $R_1$ согласованно со структурой представления на множестве $X$, то мы можем построить продолжение этого отображения, которое является эндоморфизмом представления $f$.

**Определение 3.13.** Пусть $X$ - множество образующих представления $f$. Пусть $R$ - эндоморфизм представления $f$. Множество координат $W(f, X, R \circ X)$ называется **координатами эндоморфизма представления**. $\qquad\square$

---

[4]В теореме 4.7, требования к множеству образующих более жёсткие. Поэтому теорема 4.7 говорит о продолжении произвольного отображения. Более подробный анализ дан в замечании 4.9.

**Определение 3.14.** Пусть $X$ - множество образующих представления $f$. Пусть $R$ - эндоморфизм представления $f$. Пусть $m \in M$. Мы определим **суперпозицию координат** эндоморфизма представления $f$ и элемента $m$ как координаты, определённые согласно правилу

$$(3.13) \qquad W(f, X, m) \circ W(f, X, R \circ X) = W(f, X, R \circ m)$$

Пусть $Y \subset M$. Мы определим суперпозицию координат эндоморфизма представления $f$ и множества $Y$ согласно правилу

$$(3.14) \qquad W(f, X, Y) \circ W(f, X, R \circ X) = (W(f, X, m) \circ W(f, X, R \circ X), m \in Y)$$

$\square$

**Теорема 3.15.** *Эндоморфизм $R$ представления*

$$f : A \to {}^*M$$

*порождает отображение координат представления*

$$(3.15) \qquad W[f, X, R] : W(f, X) \to W(f, X)$$

*такое, что*

$$(3.16) \qquad \begin{aligned} W(f, X, m) &\to W[f, X, R] \star W(f, X, m) = W(f, X, R \circ m) \\ W[f, X, R] &\star W(f, X, m) = W(f, X, m) \circ W(f, X, R \circ X) \end{aligned}$$

*Доказательство.* Согласно замечанию [3.8](#), мы можем рассматривать равенства ([3.13](#)), ([3.15](#)) относительно заданного множества образующих $X$. При этом координатам $W(f, X, m)$ соответствует слово

$$(3.17) \qquad W(f, X, m) \circ X = w(f, X, m)$$

а координатам $W(f, X, R \circ m)$ соответствует слово

$$(3.18) \qquad W(f, X, R \circ m) \circ X = w(f, X, R \circ m)$$

Поэтому для того, чтобы доказать теорему, нам достаточно показать, что отображению $W[f, X, R]$ соответствует отображение $w[f, X, R]$. Мы будем доказывать это утверждение индукцией по сложности $\Omega_2$-слова.

Если $m \in X$, $m' = R \circ m$, то, согласно равенствам ([3.17](#)), ([3.18](#)), отображения $W[f, X, R]$ и $w[f, X, R]$ согласованы.

Пусть для $m_1, ..., m_n \in X$ отображения $W[f, X, R]$ и $w[f, X, R]$ согласованы. Пусть $\omega \in \Omega_2(n)$. Согласно теореме [3.9](#)

$$(3.19) \qquad W(f, X, m_1...m_n\omega) = W(f, X, m_1)...W(f, X, m_n)\omega$$

Так как $R$ - эндоморфизм $\Omega_2$-алгебры $M$, то из равенства ([3.19](#)) следует

$$(3.20) \qquad \begin{aligned} W(f, X, R \circ (m_1...m_n\omega)) &= W(f, X, (R \circ m_1)...(R \circ m_n)\omega) \\ &= W(f, X, R \circ m_1)...W(f, X, R \circ m_n)\omega \end{aligned}$$

Из равенств ([3.19](#)), ([3.20](#)) и предположения индукции следует, что отображения $W[f, X, R]$ и $w[f, X, R]$ согласованы для $m = m_1...m_n\omega$.

Пусть для $m_1 \in M$ отображения $W[f, X, R]$ и $w[f, X, R]$ согласованы. Пусть $a \in A$. Согласно теореме [3.10](#)

$$(3.21) \qquad W(f, X, am_1) = aW(f, X, m_1)$$

Так как $R$ - эндоморфизм представления $f$, то из равенства (3.21) следует

$$(3.22) \qquad W(f, X, R \circ (am_1)) = W(f, X, a\,R \circ m_1) = aW(f, X, R \circ m_1)$$

Из равенств (3.21), (3.22) и предположения индукции следует, что отображения $W[f, X, R]$ и $w[f, X, R]$ согласованы для $m = am_1$.

**Следствие 3.16.** *Пусть $X$ - множество образующих представления $f$. Пусть $R$ - эндоморфизм представления $f$. Отображение $W[f, X, R]$ является эндоморфизмом представления $\Omega_1$-алгебры $A$ в $\Omega_2$-алгебре $W(f, X)$.* $\qquad\square$

В дальнейшем мы будем отождествлять отображение $W[f, X, R]$ и множество координат $W(f, X, R \circ X)$.

**Теорема 3.17.** *Пусть $X$ - множество образующих представления $f$. Пусть $R$ - эндоморфизм представления $f$. Пусть $Y \subset M$. Тогда*

$$(3.23) \qquad W(f, X, Y) \circ W(f, X, R \circ X) = W(f, X, R \circ Y)$$

$$(3.24) \qquad W[f, X, R] \star W(f, X, Y) = W(f, X, R \circ Y)$$

*Доказательство.* Равенство (3.23) является следствием равенства

$$R \circ Y = (R \circ m, m \in Y)$$

а также равенств (3.13), (3.14). Равенство (3.24) является следствием равенств (3.23), (3.16). $\qquad\square$

**Теорема 3.18.** *Пусть $X$ - множество образующих представления $f$. Пусть $R$, $S$ - эндоморфизмы представления $f$. Тогда*

$$(3.25) \qquad W(f, X, S \circ X) \circ W(f, X, R \circ X) = W(f, X, R \circ S \circ X)$$

$$(3.26) \qquad W[f, X, R] \star W[f, X, S] = W[f, X, R \circ S]$$

*Доказательство.* Равенство (3.25) следует из равенства (3.23), если положить $Y = S \circ X$. Равенство (3.26) следует из равенства (3.25) и цепочки равенств

$$
\begin{aligned}
&(W[f, X, R] \star W[f, X, S]) \star W(f, X, Y) \\
&= W[f, X, R] \star (W[f, X, S] \star W(f, X, Y)) \\
&= (W(f, X, Y) \circ W(f, X, S \circ X)) \circ W(f, X, R \circ X) \\
&= W(f, X, Y) \circ (W(f, X, S \circ X) \circ W(f, X, R \circ X)) \\
&= W(f, X, Y) \circ W(f, X, R \circ S \circ X) \\
&= W(f, X, R \circ S) \star W(f, X, Y)
\end{aligned}
$$

$\qquad\square$

Концепция суперпозиции координат очень проста и напоминает своеобразную машину Тюринга. Если элемент $m \in M$ имеет вид

$$m = m_1...m_n\omega$$

или

$$m = am_1$$

то мы ищем координаты элементов $m_i$, для того чтобы подставить их в соответствующее выражение. Как только элемент $m \in M$ принадлежит множеству образующих $\Omega_2$-алгебры $M$, мы выбираем координаты соответствующего элемента из второго множителя. Поэтому мы требуем, чтобы второй множитель в суперпозиции был множеством координат образа множества образующих $X$.

Мы можем обобщить определение суперпозиции координат и предположить, что один из множителей является множеством $\Omega_2$-слов. Соответственно, определение суперпозиции координат имеет вид

$$w(f, X, Y) \circ W(f, X, R \circ X) = W(f, X, Y) \circ w(f, X, R \circ X) = w(f, X, R \circ Y)$$

Следующие формы записи образа множества $Y$ при эндоморфизме $R$ эквивалентны.

$$(3.27) \qquad R \circ Y = W(f, X, Y) \circ (R \circ X) = W(f, X, Y) \circ (W(f, X, R \circ X) \circ X)$$

Из равенств (3.23), (3.27) следует, что

$$(3.28) \quad (W(f, X, Y) \circ W(f, X, R \circ X)) \circ X = W(f, X, Y) \circ (W(f, X, R \circ X) \circ X)$$

Равенство (3.28) является законом ассоциативности для операции композиции и позволяет записать выражение

$$W(f, X, Y) \circ W(f, X, R \circ X) \circ X$$

без использования скобок.

Рассмотрим равенство (3.25) где мы видим изменение порядка эндоморфизмов в суперпозиции координат. Это равенство также следует из цепочки равенств, где мы можем непосредственно видеть, когда изменяется порядок эндоморфизмов

$$(3.29)
\begin{aligned}
&W(f, X, m) \circ W(f, X, R \circ S \circ X) \circ X \\
&= (R \circ S) \circ m = R \circ (S \circ m) \\
&= R \circ ((W(f, X, m) \circ W(f, X, S \circ X)) \circ X) \\
&= W(f, X, m) \circ W(f, X, S \circ X) \circ W(f, X, R \circ X) \circ X
\end{aligned}$$

Из равенства (3.29) следует, что координаты эндоморфизма действуют на координаты элемента $\Omega_2$-алгебры $M$ справа.

**Определение 3.19.** Пусть $X \subset M$ - множество образующих представления

$$f : A \to {}^* M$$

Пусть отображение

$$H : M \to M$$

является эндоморфизмом представления $f$. Пусть множество $X' = H \circ X$ является образом множества $X$ при отображении $H$. Эндоморфизм $H$ представления $f$ называется **невырожденным на множестве образующих** $X$, если множество $X'$ является множеством образующих представления $f$. В противном случае, эндоморфизм $H$ представления $f$ называется **вырожденным на множестве образующих** $X$, $\qquad\square$

**Определение 3.20.** Эндоморфизм представления $f$ называется **невырожденным**, если он невырожден на любом множестве образующих. $\qquad\square$

**Теорема 3.21.** *Автоморфизм $R$ представления*

$$f : A \to {}^*M$$

*является невырожденным эндоморфизмом.*

*Доказательство.* Пусть $X$ - множество образующих представления $f$. Пусть $X' = R \circ X$.

Согласно теореме 3.7 эндоморфизм $R$ порождает отображение $\Omega_2$-слов $w[f, X, R]$.

Пусть $m' \in M$. Так как $R$ - автоморфизм, то существует $m \in M$, $R \circ m = m'$. Согласно определению 3.6, $w(f, X, m)$ - $\Omega_2$-слово, представляющее $m$ относительно множества образующих $X$. Согласно теореме 3.7, $w(f, X', m')$ - $\Omega_2$-слово, представляющее $m'$ относительно множества образующих $X'$

$$w(f, X', m') = w[f, X, R](w(f, X, m))$$

Следовательно, $X'$ - множество образующих представления $f$. Согласно определению 3.20, автоморфизм $R$ - невырожден. $\square$

## 4. Базис представления

**Определение 4.1.** Если множество $X \subset M$ является множеством образующих представления $f$, то любое множество $Y$, $X \subset Y \subset M$ также является множеством образующих представления $f$. Если существует минимальное множество $X$, порождающее представление $f$, то такое множество $X$ называется **базисом представления** $f$. $\square$

**Теорема 4.2.** *Множество образующих $X$ представления $f$ является базисом тогда и только тогда, когда для любого $m \in X$ множество $X \setminus \{m\}$ не является множеством образующих представления $f$.*

*Доказательство.* Пусть $X$ - множество образующих расслоения $f$. Допустим для некоторого $m \in X$ существует $\Omega_2$-слово

$$(4.1) \qquad w = w(f, X \setminus \{m\}, m)$$

Рассмотрим элемент $m'$, для которого $\Omega_2$-слово

$$(4.2) \qquad w' = w(f, X, m')$$

зависит от $m$. Согласно определению 3.6, любое вхождение $m$ в $\Omega_2$-слово $w'$ может быть заменено $\Omega_2$-словом $w$. Следовательно, $\Omega_2$-слово $w'$ не зависит от $m$, а множество $X \setminus \{m\}$ является множеством образующих представления $f$. Следовательно, $X$ не является базисом расслоения $f$. $\square$

*Замечание 4.3.* Доказательство теоремы 4.2 даёт нам эффективный метод построения базиса представления $f$. Выбрав произвольное множество образующих, мы шаг за шагом исключаем те элементы множества, которые имеют координаты относительно остальных элементов множества. Если множество образующих представления бесконечно, то рассмотренная операция может не иметь последнего шага. Если представление имеет конечное множество образующих, то за конечное число шагов мы можем построить базис этого представления.

Как отметил Кон в [4], стр. 96, 97, представление может иметь неэквивалентные базисы. Например, циклическая группа шестого порядка имеет базисы $\{a\}$ и $\{a^2, a^3\}$, которые нельзя отобразить один в другой эндоморфизмом представления. $\square$

*Замечание* 4.4. Мы будем записывать базис также в виде

$$X = (x, x \in X)$$

Если базис - конечный, то мы будем также пользоваться записью

$$X = (x_i, i \in I) = (x_1, ..., x_n)$$

$\square$

*Замечание* 4.5. Мы ввели $\Omega_2$-слово элемента $x \in M$ относительно множества образующих $X$ в определении 3.6. Из теоремы 4.2 следует, что если множество образующих $X$ не является базисом, то выбор $\Omega_2$-слова относительно множества образующих $X$ неоднозначен. Но даже если множество образующих $X$ является базисом, то представление $m \in M$ в виде $\Omega_2$-слова неоднозначно. Если $m_1, ..., m_n$ - $\Omega_2$-слова, $\omega \in \Omega_2(n)$ и $a \in A$, то[5]

$$(4.3) \qquad a(m_1...m_n\omega) = (am_1)...(am_n)\omega$$

Возможны равенства, связанные со спецификой представления. Например, если $\omega$ является операцией $\Omega_1$-алгебры $A$ и операцией $\Omega_2$-алгебры $M$, то мы можем потребовать, что $\Omega_2$-слова $a_1...a_n\omega x$ и $a_1x...a_nx\omega$ описывают один и тот же элемент $\Omega_2$-алгебры $M$.[6]

Помимо перечисленных равенств в $\Omega_2$-алгебре могут существовать соотношения вида

$$(4.4) \qquad w_1(f, X, m) = w_2(f, X, m) \quad m \notin X$$

Особенностью равенства (4.4) является то, что это равенство не может быть приведенно к более простому виду.[7]

На множестве $\Omega_2$-слов $w(f, X)$, перечисленные выше равенства определяют **отношение эквивалентности** $\rho(f)$, **порождённое представлением** $f$. Очевидно, что для произвольного $m \in M$ выбор соответствующего $\Omega_2$-слова однозначен с точностью до отношения эквивалентности $\rho(f)$. Тем не менее, если в процессе построений мы получим равенство двух $\Omega_2$-слов относительно заданного базиса, мы можем утверждать, что эти $\Omega_2$-слова равны, не заботясь об отношении эквивалентности $\rho(f)$.

Аналогичное замечание касается отображения $W(f, X, m)$, определённого в замечании 3.8.[8]

$\square$

**Теорема 4.6.** *Автоморфизм представления $f$ отображает базис представления $f$ в базис.*

---

[5]Например, пусть $\{\overline{e}_1, \overline{e}_2\}$ - базис векторного пространства над полем $k$. Равенство (4.3) принимает форму закона дистрибутивности

$$a(b^1\overline{e}_1 + b^2\overline{e}_2) = (ab^1)\overline{e}_1 + (ab^2)\overline{e}_2$$

[6]Для векторного пространства это требование принимает форму закона дистрибутивности

$$(a + b)\overline{e}_1 = a\overline{e}_1 + b\overline{e}_1$$

[7]Смотри, например, подразделы [3]-5.5.2, [3]-5.5.3.

[8]Если базис векторного пространства - конечен, то мы можем представить базис в виде матрицы строки

$$\overline{\overline{e}} = \begin{pmatrix} \overline{e}_1 & ... & \overline{e}_2 \end{pmatrix}$$

*Доказательство.* Пусть отображение $R$ - автоморфизм представления $f$. Пусть множество $X$ - базис представления $f$. Пусть $X' = R \circ X$.

Допустим множество $X'$ не является базисом. Согласно теореме [4.2](#) существует $m' \in X'$ такое, что $X' \backslash \{x'\}$ является множеством образующих представления $f$. Согласно теореме [3]-[2.3.5](#) отображение $R^{-1}$ является автоморфизмом представления $f$. Согласно теореме [3.21](#) и определению [3.20](#), множество $X \backslash \{m\}$ является множеством образующих представления $f$. Полученное противоречие доказывает теорему. □

**Теорема 4.7.** *Пусть $X$ - базис представления $f$. Пусть*

$$R_1 : X \to X'$$

*произвольное отображение множества $X$. Рассмотрим отображение $\Omega_2$-слов*

$$w[f, X, R_1] : w(f, X) \to w(f, X')$$

*удовлетворяющее условиям* (1), (2), (3) *теоремы* [3.7](#), *и такое, что*

$$x \in X => w[f, X, R_1](x) = R_1(x)$$

*Существует единственный эндоморфизм представления $f$[9]*

$$R : M \to M$$

*определённый правилом*

$$R \circ m = w[f, X, R_1](w(f, X, m))$$

*Доказательство.* Утверждение теоремы является следствием теорем [3.7](#), [3.11](#). □

**Следствие 4.8.** *Пусть $X$, $X'$ - базисы представления $f$. Пусть $R$ - автоморфизм представления $f$ такой, что $X' = R \circ X$. Автоморфизм $R$ определён однозначно.* □

*Замечание* 4.9. Теорема [4.7](#), так же как и теорема [3.11](#), является теоремой о продолжении отображения. Однако здесь $X$ - не произвольное множество образующих представления, а базис. Согласно замечанию [4.3](#), мы не можем определить координаты любого элемента базиса через остальные элементы этого же базиса. Поэтому отпадает необходимость в согласованности отображения базиса с представлением. □

---

Мы можем представить отображение $W(f, \overline{\overline{e}}, \overline{v})$ в виде матрицы столбца

$$W(f, \overline{\overline{e}}, \overline{v}) = \begin{pmatrix} v^1 \\ ... \\ v^n \end{pmatrix}$$

Тогда

$$W(f, \overline{\overline{e}}, \overline{v})(\overline{\overline{e}}') = W(f, \overline{\overline{e}}, \overline{v}) \begin{pmatrix} \overline{e}'_1 & ... & \overline{e}'_n \end{pmatrix} = \begin{pmatrix} v^1 \\ ... \\ v^n \end{pmatrix} \begin{pmatrix} \overline{e}'_1 & ... & \overline{e}'_n \end{pmatrix}$$

имеет вид произведения матриц.

[9]Это утверждение похоже на теорему [1]-1, с. 104.

**Теорема 4.10.** *Набор координат $W(f, X, X)$ соответствует тождественному преобразованию*

$$W[f, X, E] = W(f, X, X)$$

*Доказательство.* Утверждение теоремы следует из равенства

$$m = W(f, X, m) \circ X = W(f, X, m) \circ W(f, X, X) \circ X$$

$\square$

**Теорема 4.11.** *Пусть $W(f, X, R \circ X)$ - множество координат автоморфизма $R$. Определено множество координат $W(f, R \circ X, X)$, соответствующее автоморфизму $R^{-1}$. Множество координаты $W(f, R \circ X, X)$ удовлетворяют равенству*

$$(4.5) \qquad W(f, X, R \circ X) \circ W(f, R \circ X, X) = W(f, X, X)$$

$$W[f, X, R^{-1}] = W[f, X, R]^{-1} = W(f, R \circ X, X)$$

*Доказательство.* Поскольку $R$ - автоморфизм представления $f$, то, согласно теореме 4.6, множество $R \circ X$ - базис представления $f$. Следовательно, определено множество координат $W(f, R \circ X, X)$. Равенство (4.5) следует из цепочки равенств

$$W(f, X, R \circ X) \circ W(f, R \circ X, X) = W(f, X, R \circ X) \circ W(f, X, R^{-1} \circ X)$$
$$= W(f, X, R^{-1} \circ R \circ X) = W(f, X, X)$$

$\square$

*Замечание* 4.12. В $\Omega_2$-алгебре $M$ не существует универсального алгоритма определения множества координат $W(f, R \circ X, X)$ для заданного множества $W(f, X, R \circ X)$.[10] Мы полагаем, что в теореме 4.11 этот алгоритм задан неявно. Очевидно также, что множество $\Omega_2$-слов

$$(4.6) \qquad W(f, X, R \circ X) \circ W(f, R \circ X, X) \circ X$$

вообще говоря, не совпадает с множеством $\Omega_2$-слов

$$(4.7) \qquad W(f, X, X) \circ X$$

Теорема 4.11 утверждает, что множества $\Omega_2$-слов (4.6) и (4.7) совпадают с точностью до отношения эквивалентности, порождённой представлением $f$. $\square$

**Теорема 4.13.** *Группа автоморфизмов $G(f)$ эффективного представления $f$ в $\Omega_2$-алгебре $M$ порождает эффективное представление в $\Omega_2$-алгебре $M$.*

*Доказательство.* Из следствия 4.8 следует, что если автоморфизм $R$ отображает базис $X$ в базис $X'$, то множество координат $W(f, X, X')$ однозначно определяет автоморфизм $R$. Из теоремы 3.15 следует, что множество координат $W(f, X, X')$ определяет правило отображения координат относительно базиса $X$ при автоморфизме представления $f$. Из равенства (3.27) следует, что автоморфизм $R$ действует справа на элементы $\Omega_2$-алгебры $M$. Из равенства (3.25) следует, что представление группы ковариантно справа. Согласно теореме 4.10 набор координат $W(f, X, X)$ соответствует тождественному преобразованию. Из теоремы 4.11 следует, что набор координат $W(f, R \circ X, X)$ соответствует преобразованию, обратному преобразованию $W(f, X, R \circ X)$. $\square$

---

[10]В векторном пространстве линейному преобразованию соответствует матрица чисел. Соответственно, обратному преобразованию соответствует обратная матрица.

## 5. Многообразие базисов представления

Множество $\mathcal{B}(f)$ базисов представления $f$ называется **многообразием базисов представления** $f$.

Согласно теореме 4.6, автоморфизм $R$ представления $f$ порождает преобразование

$$(5.1) \qquad R : Y \to R \circ Y$$
$$W(f, X, R \circ Y) = W(f, X, Y) \circ W(f, X, R \circ X)$$

многообразия базисов представления. Это преобразование называется **активным**. Согласно теореме [3]-2.3.5, определено **активное представление**

$$A(f) : G(f) \to \mathcal{B}^*(f)$$

группы $G(f)$ в многообразии базисов $\mathcal{B}(f)$. Согласно следствию 4.8, это представление однотранзитивно.

*Замечание* 5.1. Согласно замечанию 4.3, могут существовать базисы представления $f$, не связанные активным преобразованием. В этом случае мы в качестве многообразия базисов будем рассматривать орбиту выбранного базиса. Следовательно, представление $f$ может иметь различные многообразия базисов. Мы будем предполагать, что мы выбрали многообразие базисов. $\qquad\square$

**Определение 5.2.** Автоморфизм $S$ представления $A(f)$ называется **пассивным преобразованием многообразия базисов представления**. Мы будем пользоваться записью

$$S(X) = S \star X$$

для обозначения образа базиса $X$ при пассивном преобразовании $S$. $\qquad\square$

**Теорема 5.3.** *Пассивное преобразование базиса имеет вид*

$$(5.2) \qquad S : Y \to S \star Y$$
$$W(f, X, S \star Y) = W(f, X, S \star X) \circ W(f, X, Y)$$

*Доказательство.* Согласно равенству (5.1), активное преобразование действует на координаты базиса справа. Равенство (5.2) следует из теорем [3]-2.5.12, [3]-2.5.13, [3]-2.5.15, согласно которым пассивное преобразование действует на координаты базиса слева. $\qquad\square$

**Теорема 5.4.** *Пусть $S$ - пассивное преобразование многообразия базисов представления $f$. Пусть $X$ - базис представления $f$, $X' = S \star X$. Пусть для базиса $Y$ существует активное преобразование $R$ такое, что $Y = R \circ X$. Положим $Y' = R \circ X'$. Тогда $S \star Y = Y'$.*

*Доказательство.* Согласно равенству (5.1), активное преобразование координат базиса $Y$ имеет вид

$$(5.3) \qquad R : X' = S \star X \to Y'$$
$$W(f, X, Y') = W(f, X, X') \circ W(f, X, Y)$$

Пусть $Y'' = S \star Y$. Из равенства (5.2) следует, что

$$(5.4) \qquad S : Y = R \circ X \to Y''$$
$$W(f, X, Y'') = W(f, X, X') \circ W(f, X, Y)$$

Из совпадения выражений в равенствах (5.3), (5.4) следует, что $Y' = Y'' = S \star Y$. Следовательно, коммутативна диаграмма

$\square$

**Теорема 5.5.** *Существует однотранзитивное* **пассивное представление**

$$P(f) : G(f) \to {}^*\mathcal{B}(f)$$

*Доказательство.* Поскольку $A(f)$ - однотранзитивное представление группы $G(f)$, то согласно определению 5.2 и теореме [3]-2.3.5, множество пассивных преобразований порождает группу. $\square$

Пусть пассивное преобразование $S$ отображает базис $Y$ в базис

$$(5.5) \qquad Y' = S \star Y$$
$$W(f, X, Y') = W(f, X, S \star X) \circ W(f, X, Y)$$

Пусть пассивное преобразование $T$ отображает базис $Y'$ в базис

$$(5.6) \qquad Y'' = T \star Y'$$
$$W(f, X, Y'') = W(f, X, T \star X) \circ W(f, X, Y')$$

Обозначим $T \star S$ произведение пассивных преобразований $T$ и $S$

$$(5.7) \qquad Y'' = T \star S \star Y$$
$$W(f, X, Y'') = W(f, X, T \star S \star X) \circ W(f, X, Y)$$

Из равенств (5.5), (5.6) следует, что

$$(5.8) \qquad Y'' = T \star S \star Y$$
$$W(f, X, Y'') = W(f, X, T \star X) \circ W(f, X, S \star X) \circ W(f, X, Y)$$

Из равенств (5.7), (5.8) и теоремы 5.5 следует, что

$$(5.9) \qquad W(f, X, T \star S \star X) = W(f, X, T \star X) \circ W(f, X, S \star X)$$

## 6. Геометрический объект представления универсальной алгебры

Активное преобразование изменяет базисы и элементы $\Omega_2$-алгебры согласовано и координаты элемента относительно базиса не меняются. Пассивное преобразование меняет только базис, и это ведёт к изменению координат элемента относительно базиса.

Допустим пассивное преобразование $S \in G(f)$ отображает базис $Y \in \mathcal{B}(f)$ в базис $Y' \in \mathcal{B}(f)$

$$(6.1) \qquad W(f, X, Y') = W(f, X, S \star X) \circ W(f, X, Y)$$

Допустим элемент $m \in M$ имеет $\Omega_2$-слово

$$(6.2) \qquad m = W(f, Y, m) \circ Y$$

относительно базиса $Y$ и имеет $\Omega_2$-слово

$$(6.3) \qquad m = W(f, Y', m) \circ Y'$$

относительно базиса $Y'$. Из (6.1) и (6.3) следует, что

$$(6.4) \qquad \begin{aligned} W(f, Y, m) \circ Y &= W(f, Y', m) \circ Y' = W(f, Y', m) \circ W(f, Y, Y') \circ Y \\ &= W(f, Y', m) \circ W(f, Y, S \star Y) \circ Y \end{aligned}$$

Сравнивая (6.2) и (6.4) получаем, что

$$(6.5) \qquad W(f, Y, m) = W(f, Y', m) \circ W(f, Y, S \star Y)$$

Так как $S$ - автоморфизм представления, то из (6.5) и теоремы 4.11 следует

$$(6.6) \qquad W(f, Y', m) = W(f, Y, m) \circ W(f, S \star Y, Y)$$

Преобразование координат (6.6) не зависит от элемента $m$ или базиса $Y$, а определённо исключительно координатами элемента $m$ относительно базиса $Y$.

**Теорема 6.1.** *Преобразования координат* (6.6) *порождают контравариант­ное эффективное $T\star$-представление группы $G(f)$, называемое* **координатным представлением в $\Omega_2$-алгебре**.

*Доказательство.* Согласно следствию 3.16, преобразование (6.6) является эн­доморфизмом представления $\Omega_1$-алгебры $A$ в $\Omega_2$-алгебре $W(f, X)$.[11]

Допустим мы имеем два последовательных пассивных преобразования $S$ и $T$. Преобразование координат (6.6) соответствует пассивному преобразованию $S$. Преобразование координат

$$(6.7) \qquad W(f, Y'', m) = W(f, Y', m) \circ W(f, T \star Y, Y)$$

соответствует пассивному преобразованию $T$. Произведение преобразований координат (6.6) и (6.7) имеет вид

$$(6.8) \qquad \begin{aligned} W(f, Y'', m) &= W(f, Y, m) \circ W(f, S \star Y, Y) \circ W(f, T \star Y, Y) \\ &= W(f, Y, m) \circ W(f, T \star S \star Y, Y) \end{aligned}$$

---

[11] Это преобразование не порождает эндоморфизма представления $f$. Координаты меня­ются, поскольку меняется базис, относительно которого мы определяем координаты. Однако элемент $\Omega_2$-алгебры $M$, координаты которого мы рассматриваем, не меняется.

и является координатным преобразованием, соответствующим пассивному преобразованию $T \star S$. Это доказывает, что преобразования координат порождают контравариантное $T \star$-представление группы $G(f)$.

Если координатное преобразование не изменяет координаты выбранного базиса, то ему соответствует единица группы $G(f)$, так как пассивное представление однотранзитивно. Следовательно, координатное представление эффективно. $\qquad\square$

Пусть $f$ - представление $\Omega_1$-алгебры $A$ в $\Omega_2$-алгебре $M$. Пусть $g$ - представление $\Omega_1$-алгебры $A$ в $\Omega_3$-алгебре $N$. Пассивное представление $g_p$ согласовано с пассивным представлением $f_p$, если существует гомоморфизм $h$ группы $G(f)$ в группу $G(g)$. Рассмотрим диаграмму[12]

$$
\begin{array}{ccc}
{}^*\mathcal{B}(f) & \xrightarrow{\quad H \quad} & {}^*\mathcal{B}(g) \\
\uparrow{\scriptstyle P(f)} & \nearrow{\scriptstyle f'} & \uparrow{\scriptstyle P(g)} \\
G(f) & \xrightarrow{\quad h \quad} & G(g)
\end{array}
$$

Так как отображения $P(f)$, $P(g)$ являются изоморфизмами группы, то отображение $H$ является гомоморфизмом групп. Следовательно, отображение $f'$ является представлением группы $G(f)$ в многообразии базисов $\mathcal{B}(g)$. Согласно построению, пассивному преобразованию $S$ $\Omega_2$-алгебры $M$ соответствует пассивное преобразование $H(S)$ $\Omega_3$-алгебры $N$.

(6.9)
$$Y'_N = H(S) \star Y_N$$

Тогда координатное преобразование в $N$ принимает вид

(6.10)
$$W(g, Y'_P, m) = W(g, Y_N, m) \circ W(g, H(S) \star Y_P, Y_P)$$

**Определение 6.2.** Мы будем называть орбиту

$$(W(g, Y_N, n) \circ W(g, H(S) \star Y_N, Y_N), S \star Y_M, S \in G(f))$$

**геометрическим объектом в координатном представлении, определённом в $\Omega_2$-алгебре** $M$. Для любого базиса $Y'_M = S \star Y_M$ соответствующая точка (6.10) орбиты определяет **координаты геометрического объекта в координатном пространстве представления** относительно базиса $Y'_M$. $\qquad\square$

**Определение 6.3.** Мы будем называть орбиту

$$(W(g, Y_N, m) \circ W(g, H(S) \star Y_N, Y_N), H(S) \star Y_N, S \star Y_M, S \in G(f))$$

**геометрическим объектом, определённым в $\Omega_2$-алгебре** $M$. Для любого базиса $Y'_M = S \star Y_M$ соответствующая точка (6.10) орбиты определяет **координаты геометрического объекта в $\Omega_2$-алгбре** $M$ относительно базиса $Y'_M$ и соответствующий вектор

$$p = W(g, Y'_P, p) \circ Y'_P$$

называется **представителем геометрического объекта в $\Omega_2$-алгбре** $M$.
$\qquad\square$

---

[12]Мы можем ослабить определение согласованности представлений и предположить, что $N$ - $\Omega_2$-алгебра. Тогда нам достаточно потребовать, чтобы отображение $(h, H)$ было морфизмом представлений. Но тогда элементы аффинного представления не могут быть геометрическим объектом векторного пространства.

Так как геометрический объект - это орбита представления, то согласно теореме [3]-2.4.13 определение геометрического объекта корректно.

Мы будем также говорить, что $p$ - это **геометрический объект типа** $H$

Определение 6.2 строит геометрический объект в координатном пространстве. Определение 6.3 предполагает, что мы выбрали базис представления $g$. Это позволяет использовать представитель геометрического объекта вместо его координат.

**Теорема 6.4 (принцип инвариантности).** *Представитель геометрического объекта не зависит от выбора базиса $Y'_M$.*

*Доказательство.* Чтобы определить представителя геометрического объекта, мы должны выбрать базис $Y_M$, базис $Y_P$ и координаты геометрического объекта $W(g, Y_P, p)$. Соответствующий представитель геометрического объекта имеет вид

$$p = W(g, Y_P, p) \circ Y_P$$

Базис $Y'_M$ связан с базисом $Y_M$ пассивным преобразованием

$$Y'_M = S \circ Y_M$$

Согласно построению это порождает пассивное преобразование (6.9) и координатное преобразование (6.10). Соответствующий представитель геометрического объекта имеет вид

$$\begin{aligned}
p' &= W(g, Y'_P, p') \circ Y'_P \\
&= W(g, Y_P, p) \circ W(g, H(S) \star Y_P, Y_P) \circ W(g, Y_P, H(S) \star Y_P) \circ Y_P \\
&= W(g, Y_P, p) \circ Y_P = p
\end{aligned}$$

Следовательно, представитель геометрического объекта инвариантен относительно выбора базиса. $\square$

**Теорема 6.5.** *Множество геометрических объектов типа $H$ является $\Omega_3$-алгеброй.*

*Доказательство.* Пусть

$$p_i = W(g, Y_P, p_i) \circ Y_P \quad i = 1, ..., n$$

Для операции $\omega \in \Omega_3(n)$ мы положим

$$(6.11) \qquad p_1...p_n\omega = W(g, Y_P, p_1)...W(g, Y_P, p_n)\omega \circ Y_P$$

Так как для произвольного эндоморфизма $S$ $\Omega_2$-алгебры $M$ отображение $W(g, H(S)\star Y_N, Y_N)$ является эндоморфизмом $\Omega_3$-алгебры $N$, то определение (6.11) корректно. $\square$

**Теорема 6.6.** *Определено представление $\Omega_1$-алгебры $A$ в $\Omega_3$-алгебре $N$ геометрических объектов типа $H$.*

*Доказательство.* Пусть

$$p = W(g, Y_P, p) \circ Y_P$$

Для $a \in A$ мы положим

$$(6.12) \qquad ap = aW(g, Y_P, p_1) \circ Y_P$$

Так как для произвольного эндоморфизма $S$ $\Omega_2$-алгебры $M$ отображение $W(g, H(S)\star Y_N, Y_N)$ является эндоморфизмом представления $g$, то определение (6.12) корректно. $\quad\square$

## 7. Примеры базиса представления универсальной алгебры

7.1. **Векторное пространство.** Рассмотрим векторное пространство $\overline{V}$ над полем $F$. Если дано множество векторов $\overline{e}_1$, ..., $\overline{e}_n$, то, согласно алгоритму построения координат над векторным пространством, координаты включают такие элементы как $\overline{e}_1 + \overline{e}_2$ и $a\overline{e}_1$. Рекурсивно применяя правила, приведенные в определении 3.6, мы придём к выводу, что множество векторов $\overline{e}_1$, ..., $\overline{e}_n$, порождает множество линейных комбинаций

$$a^1\overline{e}_1 + ... + a^n\overline{e}_n$$

Согласно теореме 4.2, множество векторов $\overline{e}_1$, ..., $\overline{e}_n$, является базисом при условии, если для любого $i$, $i = 1, ..., n$, вектор $\overline{e}_i$ не является линейной комбинацией остальных векторов. Это требование равносильно требованию линейной независимости векторов.

7.2. **Представление группы на множестве.** Рассмотрим представление из примера [3]-2.6.2. Мы можем рассматривать множество $M$ как объединение орбит представления группы $G$. В качестве базиса представления можно выбрать множество точек таким образом, что одна и только одна точка принадлежит каждой орбите представления. Если $X$ - базис представления, $A \in X$, $g \in G$, то $\Omega_2$-слово имеет вид $A + g$. Поскольку на множестве $M$ не определены операции, то не существует $\Omega_2$-слово, содержащее различные элементы базиса. Если представление группы $G$ однотранзитивно, то базис представления состоит из одной точки. Этой точкой может быть любая точка множества $M$.

## 8. Список литературы

## 10. Специальные символы и обозначения